\documentclass[10pt,a4paper,twoside]{amsart}

\usepackage{amsfonts, amssymb, amsmath, amsthm}
\usepackage{bm}
\usepackage{mathrsfs} 
\usepackage{enumerate}

\usepackage{url}

\usepackage{graphicx}

\newtheorem{thm}{Theorem}[section]
\newtheorem{lem}[thm]{Lemma}

\newtheorem{defi}[thm]{Definition}
\newcommand{\Z}[1]{\mathbb{Z}/#1\mathbb{Z}}

\newcommand{\RsurZ}{\mathbb{R}/\mathbb{Z}}

\let\ve=\varepsilon

\def\goes{\mathrel{\rightarrow}}
\def\order{\asymp}

\def\Bcal{\mathcal{B}}

\def\Pcal{\mathcal{P}}
\def\Xcal{\mathcal{X}}

\def\Ocal{\mathcal{O}}

\def\1{1\!\!\!1}

\def\Cscr{\mathscr{C}}

\def\mode{\mathbin{\,\textrm{mod}^*}}

\title{Cusps of primes in dense subsequences -- Bypassing the $W$-trick }
\author{Olivier Ramar\'e}

\begin{document}
\address[O. Ramar\'e]{CNRS/ Institut de Math\'ematiques de Marseille, Aix 
 Marseille Universit\'e, U.M.R. 7373, Site Sud, Campus de Luminy, Case 907, 
 13288 
 Marseille Cedex 9, France.}
 \email{olivier.ramare@univ-amu.fr}

\date{\sl December, 12th of 2024}


\maketitle

\begin{abstract}
  
  Let the $A$-cusps of a dense subset
  $\mathcal{P}^*\in[\sqrt{N},N]$ of primes be points
  $\alpha\in\RsurZ$ that are such that
  $|\sum_{\substack{p\in\mathcal{P}^*}} e(\alpha p)|\ge
  |\mathcal{P}^*|/A$. We establish that any $(1/N)$-well spaced subset
  of $A$-cusps contains at most $20A^2K\log(2A)$ points, where
  $K=N/(|\mathcal{P}^*|\log N)$.  We further show that any
  $B$-cusps~$\xi$ is accompanied, when $B\le \sqrt{A}$, by a large
  proportion of $A$-cusps of the shape $\xi+(a/q)$.  We conclude this
  study by showing that, given $A\ge2$, the characteristic function
  $\1_{\mathcal{P}^*}$ may be decomposed in the form
  $\1_{\mathcal{P}^*}=(V(z_0)\log N)^{-1}f^\flat +f^\sharp$ where the
  trigonometric polynomial of $f^\sharp$ takes only values
  $\le |\mathcal{P}^*|/A$, and~$f^\flat$ is a bounded non-negative
  function supported on the integers prime to $M$; the parameters
  $z_0$ and $M$ are given in terms of~$A$, while
  $V(z_0)=\prod_{p<z_0}(1-1/p)$. The function $f^\flat$ satisfies more
  regularity properties. In particular, its density with respect to
  the integers $\le N$ and coprime to~$M$ is again~$K$. This transfers
  questions on~$\mathcal{P}^*$ to problems on integers coprime to the
  modulus~$M$.
\end{abstract}

\section{Results in the large}

During their investigations on the  primes, B. Green in \cite[Proof of Lemma
6.1]{Green*05} and B. Green \& T. Tao in~\cite{Green-Tao*04} were led to
consider the large values of the Fourier polynomial built on some
dense subset of the primes. Let us introduce a notion for clarity.
\begin{defi}
  \label{defCusps}
  Let $\mathcal{P}^*\subset[1,N]$ be a subset of the primes, and let
  $A\ge1$ be given. We define the set of \emph{$A$-cusps} by
  \begin{equation}
    \label{eq:3}
    \Cscr(\mathcal{P}^*,A)=
    \biggl\{\alpha\in\mathbb{R}/\mathbb{Z}:
    \biggr|\sum_{\substack{p\in\mathcal{P}^*}}
    e(\alpha p)\biggl|\ge \sum_{\substack{p\in\mathcal{P}^*}}1/A
    \biggr\}.
  \end{equation}
\end{defi}
\noindent
As the involved trigonometric polynomial is continuous, the set
$\Cscr(\mathcal{P}^*, A)$ is closed, hence compact, and is more
precisely a finite union of arcs. The above set may sometimes be
called \emph{spectrum}, but, first this word is overloaded and, second
the right-hand side is often $N/A$ rather than the one we employ
above. The similarity with the same word in automorphic theory is also
welcome, see the structure theorem below.

\subsubsection*{Cusps are scarce}
Given a subset $\mathcal{P}^*\subset[\sqrt{N},N]$ of the primes, of
cardinality 
$N/(K\log N)$, B.~Green \& T.~Tao proved in \cite{Green-Tao*04} that
the cardinality of a $1/N$-well spaced subset of $\Cscr(\mathcal{P}^*,A)$ is
$\Ocal_\ve((A^2K)^{1+\ve})$, for every positive $\ve$. We sharpened
this result in~\cite{Ramare*22} to $\Ocal(A^2K\log A)$. We are
somewhat more precise in the next result.
\begin{thm}
  \label{thmCusps}
  Notation being as in Definition~\ref{defCusps}, the set
  $\Cscr(\mathcal{P}^*,A)$ is a finite union of arcs.
  Any $1/N$-well spaced subset of $\Cscr(\mathcal{P}^*,A)$ contains
  at most $(4e^\gamma+o(1))A^2K\log(2A)$ points, where
  $K=N/(|\mathcal{P}^*|\log N)$.
  When $N\ge 10^4$, such a $1/N$-well spaced set contains at
  most~$19A^2K\log(2A)$ points.
\end{thm}
\noindent
Let us recall that a set $\Xcal\subset\mathbb{R}/\mathbb{Z}$ is said
to be $\delta$-well spaced when $\min_{x\neq
  x'\in\Xcal}|x-x'|_{\mathbb{Z}}\ge\delta$, where
$|y|_{\mathbb{Z}}=\min_{k\in\mathbb{Z}}|y-k|$ denotes in a rather
unusual manner the distance to
the nearest integer.

In particular, the measure of $\Cscr(\mathcal{P}^*,A)$
is~$\Ocal(A^2(\log A)/N)$, as $N$ grows.

\subsubsection*{Conditional optimality}
\begin{thm}
  \label{L1}
  If Theorem~\ref{thmCusps} holds with the bound $\Ocal(A^2Kf(\log A))$ for
  some non-decreasing positive function $f$, then
  \begin{equation*}
    \int_0^1 \biggl|\sum_{p\in\mathcal{P}^*}e(p\alpha)\biggr|d\alpha
    \ll
    \sqrt{\frac{|\mathcal{P}^*|\,f(\log N)}{\log N}}.
  \end{equation*}
\end{thm}
\noindent
A lower bound of size $|\mathcal{P}^*|/\sqrt{N}$ can
be infered from the work~\cite{Vaughan*88} of R. C.~Vaughan. It is
reproved in a more general fashion by E.~Eckels, S.~Steven, A.~Leodan
\& R.~Tobin in \cite{Eckels-Jin-Leodan-Tobin*23}. When
restricting to the full sequence of primes, the best known upper
bound is $(\frac{\sqrt{2}}{2}+o(1))\sqrt{N/\log N}$, due to D.~A.~Goldston
in~\cite{Goldston*92b}. Both
Vaughan and Goldston tend to believe that the good order of
magnitude is~$\sqrt{N/\log N}$, in which case Theorem~\ref{thmCusps}
would be optimal. However, when~$A$ is small, we have only been able
to build examples with~$\Ocal(A^2)$ cusps (see below).

\subsubsection*{The case of the full sequence of primes}
We may approximate the trigonometric polynomial $T$ on the primes via a
local model, i.e. write
\begin{equation}
  \label{FullPol}
  T(\alpha)=\sum_{p\le N}e(p\alpha)
  =\frac{1}{V(z_0)\log N}\sum_{\substack{ n\le N\\
      (n,P(z_0))=1}}e(n\alpha)
  +\Ocal\biggl(\frac{N}{z_0\log N}\biggr)
\end{equation}
where $z_0\le \log\log N$ is a parameter at our
disposal and
\begin{equation}
  \label{eq:4}
  P(z_0)=\prod_{p< z_0}p,\quad V(z_0)=\prod_{p<z_0}\biggl(1-\frac{1}{p}\biggr).
\end{equation}
Eq.~\eqref{FullPol} is proved in Subsection~\ref{ProofFullPol}.
As a consequence and when limiting our study to~$A$ small, say
$A=o(\log\log N)$, we find that the set of~$A$-cusps of~$T$
is a union of arcs around points from $\{a/q: (a,q)=1, q|P(z_0)\}$, for
some~$z_0$ (chosen for the error term in \eqref{FullPol} to
become $<N/(A\log N)$). Around~$a/q$ and when $q|P(z_0)$, say in~$\alpha=\frac aq+\beta$, we readily find that
$|T(\alpha)|\le \min(\frac{\mu^2(q)N}{\varphi(q)\log
  N},P(z_0)/\|P(z_0)\beta\|)+\Ocal(N/(z_0\log N))$. Furthermore, on adapting the
proof of P. Bateman from \cite{Bateman*72}\footnote{The result of
  Bateman does not include the square-free condition on~$q$.}, we readily find that
\begin{equation*}
  \sum_{\substack{ \varphi(q)\ge A}}\mu^2(q)\varphi(q)
  \sim
  \frac{A^2}{2}.
\end{equation*}
Gathering our results, we find that $\Cscr(\mathcal{P}\cap[1,N],A)$ is a
union of about $(A^2/2)$~arcs of length $\Ocal(1/N)$.
\begin{figure}
  \centering
  \includegraphics[scale=0.5]{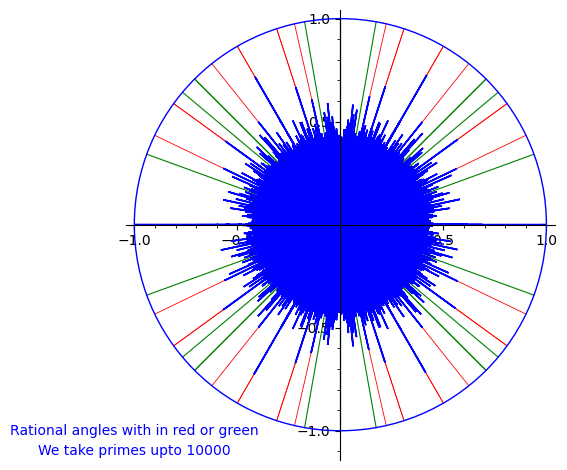}
  \caption{The modulus $|T(\alpha)|/T(0)$. In red (resp. green), the first lines
    from the origin to rational points with square-free
    (resp. non square-free) denominators.}
\end{figure}
The figure we display has been again doctored a bit, but exhibit the
presence of very sharp cusps. We produced it by using Sage, cf~\cite{sagemath}.
\subsubsection*{A subsequence exhibiting non-rational cusps}
A part of the Farey sequence, namely the points with square-free
denominators, appears as cusps  for the full sequence of primes. 
We found the presence of this sequence in many numerical
examples\footnote{We tried random subsequences of primes of relative
  density $1/2$, then Ramanujan primes and, finally we selected successively one prime
out of two.}. Let us
give an example were non-rational cusps appear. We select
\begin{equation}
  \label{eq:5}
  \mathcal{P}^*_0=\bigl\{p: \{p\sqrt{2}\}\le 1/2\bigr\}
\end{equation}
where $\{x\}$ denotes the fractional part of the real number~$x$. By
detecting the condition $\{p\sqrt{2}\}\le 1/2$ through Fourier analysis, we
readily find that
\begin{multline}
  \label{SpecialPol}
  T^*_0(\alpha)=\sum_{\substack{p\le N\\ \{p\sqrt{2}\}\le 1/2}}e(p\alpha)
  =\frac{(1/2)}{V(z_0)\log N}\sum_{\substack{2\le n\le N\\
      (n,P(z_0))=1}}e(n\alpha)
 \\ +\sum_{\substack{|h|\le H\\ \text{$h$ odd}}}\frac{1}{i\pi h}
  \frac{1}{V(z_0)\log N}\mkern-7mu\sum_{\substack{2\le n\le N\\
      (n,P(z_0))=1}}\mkern-15mu e\bigl(n{(h\sqrt{2}+\alpha)}\bigr)
  +\Ocal\biggl(\biggl(\frac{1}{H}+\frac{1}{z_0}\biggr)\frac{N}{\log N}\biggr).
\end{multline}
where $P(z_0)=\prod_{p< z_0}p$ as above. This is proved in Subsection~\ref{ProofSpecialPol}.
\begin{figure}
  \centering
  \includegraphics[scale=0.5]{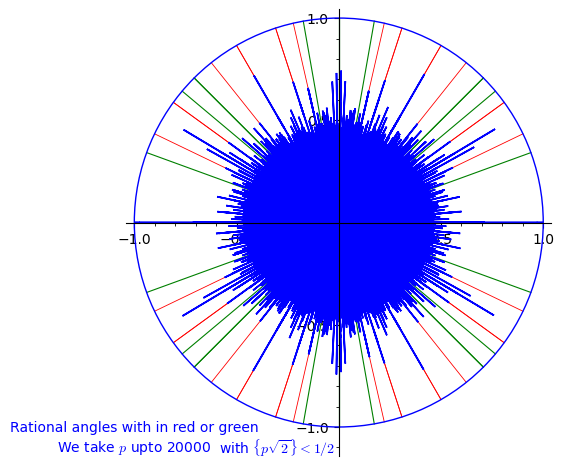}
  \caption{The modulus $|T^*_0(\alpha)|/T^*_0(0)$. In red (resp. green), the first lines
    from the origin to rational points with square-free
    (resp. non square-free) denominators: some new cusps arise.}
\end{figure}
The figure we produce has been again doctored; it also exhibits the
presence of very sharp cusps. 

\subsubsection*{Cusps and Farey points}

In the two previous examples, the Farey sequence plays a major role
which seems to be related to the coprimality condition~$(n,P(z_0))=1$
in~\eqref{FullPol} and~\eqref{SpecialPol}. The next theorem shows
that this is the general situation.
\begin{thm}
  \label{Companions}
  Notation being as in Definition~\ref{defCusps}. Let $N\ge10^4$,
  $A\in [2,\sqrt{N}]$, $B\in[1,A]$
   and $\xi\in\Cscr(\mathcal{P}^*,B)$ be given. Define $t^*(\alpha)=|T^*(\alpha)|/T^*(0)$.
  The set
  \begin{equation*}
    \mathcal{F}=\bigl\{\xi+(a/q):a\mode q, q\le A/B\bigr\}\cap \Cscr(\mathcal{P}^*,A)
  \end{equation*}
  contains more than $A^2/(6800 B^4Z^2K\log A)$ elements, where $Z\in[1/B,1]$
  is the maximum value of $t^*(\xi+(a/q))$ for $\xi+(a/q)\in\mathcal{F}$.
\end{thm}
\noindent
Lemma~\ref{iniStructureCusps} also exhibits some
easy properties on $\Cscr(\mathcal{P}^*,A)$. 

In particular, the cusps at any of 0, $1/3$, $1/2$ and $2/3$
generate $\order A^2/(K\log A)$ $A$-cusps by simply adding $a/q$
with $q\le A$ to them.

\subsubsection*{A structure theorem without $W$-trick}
The ubiquity of the Farey points with square-free denominators was
rather a surprise. We already encountered in~\cite{Ramare*22} the
major influence of the parameter $z_0$ above. Since we prove a
numerically better version in Theorem~\ref{Bel} below, the reader may
see the effect of the ``unsieving'' parameter~$z_0$. Musing on this
effect, we prove the next structure theorem, which may be seen as a
general form of~\eqref{FullPol} and~\eqref{SpecialPol}. We use the
generic notation
\begin{equation}
  \label{defSg}
  S(g,\alpha)=\sum_{n\le N}g(n)e(n\alpha).
\end{equation}
\begin{thm}
  \label{MainDecB}
  Let $\epsilon\in(0,1/3]$ and $A,K\ge1$ be given. Let
  \begin{equation}
    \label{defz0}
    z_0=\exp\bigl(25000\,A^3(\log 2A)^2K\bigr)/\epsilon.
  \end{equation}
  Let $M$ be an integer divisible by $P(z_0)$ and all whose prime factors are~$<z_0$.
  Let $N\ge (Mz_0)^{1+\epsilon^{-1}}$.
  Let $\mathcal{P}^*$ be subset of the primes inside $[\sqrt{N},N]$ of
  cardinality (at least) $N/(K\log N)$ and let $T^*(\alpha)=\sum_{p\in\mathcal{P}^*}e(p\alpha)$.
  We may decompose the characteristic function $\1_{\mathcal{P}^*}$ of
  $\mathcal{P}^*$ as
  \begin{equation}
    \label{GenPol}
    \1_{\mathcal{P}^*}=\frac{f^\flat}{V(z_0)\log N} +f^\sharp
  \end{equation}
  in such a way that the following properties hold:
  \begin{itemize}
  \item For every $\alpha\in\RsurZ$, $|S(f^\sharp,\alpha)|\le |T^*(\alpha)|$
  and $|S(f^\flat,\alpha)|\le |T^*(\alpha)|V(z_0)\log N$.
  \item For every $\alpha\in\RsurZ$, we have $|S(f^\sharp,\alpha)|<
    T^*(0)/A$.
  \item $0\le f^\flat(\ell)\le 2(1+\epsilon)^2$
    and $f^\flat(\ell)=0$ when $(\ell,M)\neq1$.
  \item For any integer $a$, we have $S(f^\flat,a/M) = T^*(a/M)V(z_0)\log N$.
  \end{itemize}
\end{thm}
\noindent
This may be compared with several results of B. Green \& T. Tao (see
also \cite{Fox-Zhao*15} by J. Fox \& Y. Zhao and
\cite{Conlon-Fox-Zhao*15} by D.Conlon, J. Fox \& Y. Zhao): the
two differences are (a) that we approximate $\1_{\mathcal{P}^*}$ and not only
a $W$-tricked version of it and, (b) that this approximation is valid
on the full unit circle. This second modification is of lesser
importance and essentially improves on the readability.
The decomposition~\eqref{GenPol} may be compared on
the Fourier side with the ones given in~\eqref{FullPol} and~\eqref{SpecialPol}. 
In short, the above decomposition now takes care of all the cusps. It
is surprising that, in spirit, taking care only of the rational cusps
leads to a full description. 
The decomposition we propose \emph{transfers} the problem on
primes from $\mathcal{P}^*$ to a problem on $P(z_0)$-sifted integers
(i.e. integers prime to $P(z_0)$) for a weighted sequence $f^\flat$.

The reader should notice that we used another normalization than the
one of \cite{Green-Tao*04}: the function $f^\flat$ on the
$P(z_0)$-sifted integers has the same density, $K$, with respect to
the $P(z_0)$-sifted integers, than $\mathcal{P}^*$ with respect of the
primes. Now the local upper bound is around~$2$, reflecting the loss
in the method, while in \cite{Green-Tao*04}, it is around~1, the loss
in the method being reflected by a lower density. This is only a
cosmetic change, but is to be noted when comparing distinct sources.
\subsubsection*{Scope and limitations}
Theorem~\ref{MainDecB} has been thought with~$A$ small, but still
allows this parameter to grow with~$N$. As of now, if $\log M\ll z_0$, we require
\begin{equation*}
  A\ll \frac{(\log\log N)^{1/3}}{(\log\log\log N)^{2/3}}.
\end{equation*}
This goes to infinity but remains very small. The exponent
$1/3$ can be reduced to $1/2$ when asking only for a decomposition
over $\Z{\Ocal(1)N}$, but this is all the present method offers.

When $A$ is close to~1 and $\Cscr(\mathcal{P}^*,A)$ contains a point away
from $1$ and $1/2$, we may consider using the methods
from~G. Fre\u \i man in~\cite{Freiman*62} or from~V. Lev in~\cite{Lev*05}.

We concentrate in this paper on the case of the initial segment of the
primes, but, as the methods are essentially combinatorial, the results
extend to any sieve setting. How so is not exactly clear as far as the
dependence in the various parameters is concerned. We delay such
generalizations to later papers, but as it is not more work, we
prepare the ground by introducing a parameter~$\tau$ in the enveloping
sieve.

\subsubsection*{Inverse question: an open problem}
Theorem~\ref{MainDecB} creates, from a non-negative function on the
primes $\1_{\mathcal{P}^*}$, a non-negative function~$f^\flat$ on the
integers sifted by~$P(z_0)$. The main question is to know whether we 
encounter only a special class of functions $f^\flat$, or whether
essentially any such function could occur. Notice that the values of~$f^\flat$ at
each integer is maybe not fully determined.

\subsubsection*{Remarks on the explicit aspect}
Following for instance B.~Green \& T.~Tao in \cite{Green-Tao*08-2}, we  specify the absolute constants. This is of no consequence and
simplifies the understanding of the statements. The reader may
forget their value in a first reading, especially so since we did not try
to optimize these. Notice however that, if these constants were
less, one could test the results against numerical trials with more
efficiency.

\subsubsection*{Acknowledgments}
A large part of this work was completed when the author was enjoying
the hospitality of the Indian Statistical Institute in Kolkata. It has
also has been partially supported by the Indo-French IRL Relax
and by the joint FWF-ANR project Arithrand: FWF: I 4945-N and ANR-20-CE91-0006.
The author also benefited from fruitful discussions with
R.~Balasubramanian.


\section{Local model expansions}

\subsection{Local model for the primes. Proof of Eq. \eqref{FullPol}}
\label{ProofFullPol}
In order to prove Eq.~\eqref{FullPol}, we show that the left and
right-hand sides coincide for every $\alpha\in\RsurZ$.
\begin{lem}
  \label{approxT}
  When $N\ge2$, $\alpha=(a/q)+\beta$ with $a$ prime to $q$ and $|q\beta|\le
  N/(\log N)^A$, we have
  \begin{equation*}
    T(\alpha)=\sum_{p\le N}e(p\alpha)
    =\frac{\mu(q)}{\varphi(q)\log N}\sum_{n\le N}e(n\beta)
    +\Ocal\biggl(\frac{N}{(\log N)^2}\biggr).
  \end{equation*}
\end{lem}

\begin{proof}
  This is a classical result, obtained for instance by combining
  \cite[Theorem 3.1]{Vaughan*81} of the book of R.C.~Vaughan together with an obvious adaptation
  of \cite[Lemma 3.1]{Vaughan*81} in the same book. Since we do not
  weight the primes by $\log p$, we have to limit our saving to
  $\Ocal(N/(\log N)^2)$.
\end{proof}

\begin{lem}
  \label{handlecoprime}
  When $\alpha=(a/q)+\beta$ with $a$ prime to $q$ and $|q\beta|\le 1/(2P(z_0))$ have
  \begin{equation*}
    \sum_{\substack{ n\le N\\
        (n,P(z_0))=1}}e(n\alpha)
    =
    \frac{\mu(q)V(z_0)}{\varphi(q)}
    \sum_{\substack{n\le N}}e(n\beta)+\Ocal\biggl(qP(z_0)+
    \frac{N|\beta|P(z_0)^2}{\varphi(q)}
    \biggr).
  \end{equation*}
\end{lem}

\begin{proof}
  Let us use the shortcut $M_0=P(z_0)$. By using the Moebius function
  to handle the coprimality condition, we readily find that
  \begin{align*}
    L(\alpha)
    &=\sum_{\substack{ n\le N\\
        (n,P(z_0))=1}}e(n\alpha)
    =\sum_{d|M_0}\mu(d) \sum_{m\le N/d}e(dm\alpha)
    \\&=
    \sum_{\substack{d|M_0\\ q|d}}\mu(d) \sum_{m\le N/d}e(dm\alpha)
    +\Ocal\biggl(\sum_{\substack{d|M_0\\ q\nmid d}} 1/\|d\alpha\|\biggr).
  \end{align*}
  In the error term, we have $1/\|d\alpha\|\ll q$. In the main term
  (which may be non-zero only when $q|M_0$),
  we may replace $\alpha$ by~$\beta$. This quantity is therefore
  independant on $a$ on which we may sum, getting
  \begin{equation*}
    L(\alpha)=\frac{\mu(q)}{\varphi(q)}L(\beta)+\Ocal(qM_0).
  \end{equation*}
  To proceed, we use the above with $q=1$ and notice that
  $e(n\beta)-e((m+k)\beta)\ll k|\beta|$. As a consequence, we find
  that
  \begin{equation*}
    d\sum_{\substack{n\le N\\ d|n}}e(n\beta)=
    \sum_{\substack{n\le N}}e(n\beta)+\Ocal(dN|\beta|)
  \end{equation*}
  from which we infer that
  $
    L(\beta)=V(z_0)\sum_{\substack{n\le
        N}}e(n\beta)+\Ocal(N|\beta|M_0^2)$, ending the proof.
  
\end{proof}

\begin{proof}[Proof of Eq.~\eqref{FullPol}]
  We compare the expressions of Lemmas~\ref{approxT}
  and~\ref{handlecoprime} to get the result.
\end{proof}

\subsection{Another local model. Proof of Eq. \eqref{SpecialPol}}
\label{ProofSpecialPol}
\begin{lem}
  \label{V1}
  Let $I\subset\RsurZ$ be an interval and $\chi_I$ be its indicator function. 
For each positive integer $H$ there exist coefficients $a_H(h)$ and
$C_h$
for $-H\le h\le H$ with $|a_H (h)| \le \min(|I|,1/(|h|\pi))$
and $|C_h|\le 1$ such that the trigonometric polynomial
\begin{equation}
\chi^*_{I,H}(t) = |I| +
\sum_{0< |h|\le H}
 a_H(h)
 e(ht)
\end{equation}
satisfies, for every $t\in\RsurZ$,
\begin{equation}
  |\chi_I(t)-\chi^*_{I,H}(t)|\le
 \frac{1}{H+1}
 \sum_{|h|\le H}
 C_h
\biggl(1-\frac{|h|}{H+1}\biggr)
 e(ht).
\end{equation}
We also have $|a_H(h)|\le \min(|I|,1/(|h|+1))$.
\end{lem}
\noindent
This is, up to a trivial change of notation, a specialisation of
~\cite[Theorem 19]{Vaaler*85} by J.~Vaaler, reproduced in~\cite[Lemma
6.2]{Madritsch-Tichy*19} by M.~Madritsch \& R.~Tichy. 

\begin{proof}[Proof of Eq.~\eqref{SpecialPol}]
  We detect the condition $\{p\sqrt{2}\}\le 1/2$ in $T_0^*(\alpha)$ by
  using Lemma~\ref{V1}. We then employ the local model
  expansion~\eqref{FullPol} of $T(\alpha)$ to infer the expression we claim.
\end{proof}

\subsection{General structure and remarks on these expansions}
\label{GenLocal}
The gist of~\eqref{FullPol} and~\eqref{SpecialPol} is, on
employing the vocabularity of algebraic number theory, to separate the
behaviour at the infinite place (size conditions) to the one at the
finite places (congruence conditions). However, and
since we are handling sequences and not points, the contribution at
the finite places do not split prime per prime. Though there are
surely limitations to this kind of expansions (see for instance
\cite{Granville-Soundararajan*01} by A.~Granville \& K.~Soundararajan), Theorem~\ref{MainDecB}
has a wide range of application.

Together with G.K. Viswanadham in \cite{Ramare-Viswanadham*24}, we
devised the same kind of expansion for integers that are sums of two squares.

\part{The initial large sieve inequality}
Most of the work in this part is devoted to proving the next
theorem, though
the enveloping sieve we build will also be used to prove the structure
theorem~\ref{MainDecB}. 
\begin{thm}
  \label{ExtensionpreciseNow}
  Let $\Xcal$ be a $\delta$-well spaced subset of
  $\mathbb{R}/\mathbb{Z}$ and $N\ge 1$.
  Let $({}u_p)_{p\le N}$ be a sequence of complex numbers.
  We have
  \begin{equation*}
    \sum_{x\in\Xcal}
    \biggl|\sum_{ p\le N}{}u_p e(xp)\biggr|^2
    \le
     (4e^\gamma+o(1)) \frac{N+\delta^{-1}}{\log N}\log(2|\Xcal|)\sum_{p\le N}|{}u_p|^2,
   \end{equation*}
   where $o(1)$
  is a function that goes to~0 when $N$ and $|\Xcal|$ go to infinity.
  This constant is at least as large as~$\frac12+o(1)$, and,
  when $N\ge 10^4$, it is at most~19.
\end{thm}
\noindent
In most applications, $\delta^{-1}$ is~$\Ocal(N)$.
A version of this bound was proved in~\cite{Ramare*22} with
$8e^\gamma$ (resp.~280) rather
than~$4e^\gamma$ (resp.~19).
The gain stems from a finer study
of the coefficients~$w_q(z;z_0)$ resulting in
Lemma~\ref{CourageousTer}. When comparing to the similar \cite[Lemma
4.2]{Ramare*22}, one sees that we save $\sqrt{q}$, resulting in the
optimal bound $(1+o(1))/z_0$.

\section{The enveloping sieve with an unsifted part}
\label{ES}
We fix two real parameters $z_0\le z$. It is easy to reproduce the
analysis of \cite[Section 3]{Ramare-Ruzsa*01} as far as exact formulae
are concerned, but one gets easily sidetracked towards slightly
different formulae. The reader may for instance compare \cite[Lemma
4.2]{Ramare*95} and \cite[(4.1.14)]{Ramare-Ruzsa*01}.  Similar
material is also the topic of \cite[Chapter 12]{Ramare*06}.  So we
present a path leading to \cite[(4.1.14)]{Ramare-Ruzsa*01} in our
special case that we extend a bit: we add a parameter~$\tau$ and
consider integers that prime to $\prod_{p<z,p\nmid\tau}p$. The case
$\tau=1$ is the one we shall require here, but such a parameter $\tau$
naturally appears when considering primes in some fixed arithmetic
progression. We may also assume that $\tau$ is prime to~$P(z_0)$.

\subsubsection*{Main players}

We define
\begin{equation}
  \label{defGdz}
  G_d(y)
  =\sum_{\substack{\ell\le y,\\
      (\ell,d)=1}}\frac{\mu^2(\ell)}{\varphi(\ell)},
  \quad
  G(y)=G_1(y).
\end{equation}
We generalize the
definition~\eqref{defGdz} by
\begin{equation}
  \label{defGdyz0}
  G_d(y;z_0)
  =\sum_{\substack{\ell\le y,\\
      (\ell,dP(z_0))=1}}\frac{\mu^2(\ell)}{\varphi(\ell)},
  \quad
  G(y;z_0)=G_1(y;z_0).
\end{equation}
We further set
\begin{equation}
  \label{defbetan}
  \beta_{z_0,z}(n;\tau)=\biggl(\sum_{\substack{d|n}}\lambda_d\biggr)^2,\qquad
  \lambda_d=
  \1_{(d,\tau P(z_0))=1}\frac{\mu(d)dG_{d\tau}(z/d;z_0)}{\varphi(d)G_\tau(z;z_0)}
  .
\end{equation}
So the reader can see that the parameter~$\tau$ is just an extension
of the~$P(z_0)$.  
Let us recall \cite[Theorem 2.1]{Ramare*22}. Section~3 of
\cite{Ramare-Ruzsa*01} corresponds to the case~$z_0=1$.
\begin{lem}
  \label{Fourierbetan}
  The coefficients $\beta_{z_0,z}(n)$ admit the expansion
  \begin{equation*}
    \beta_{z_0,z}(n,\tau)
    =\sum_{\substack{q\le z^2,\\ q| P(z)/P(z_0)\\ (q,\tau)=1}}w_q(z;z_0,\tau)c_q(n)
  \end{equation*}
  where $c_q(n)$ is the Ramanujan sum and where
  \begin{equation*}
    w_q(z;z_0,\tau)
    =
    \frac{\mu(q)}{\varphi(q)G_\tau(z;z_0)}\frac{G_{[q]}(z;z_0,\tau)}{G_\tau(z;z_0)}
  \end{equation*}
  with the definitions
  \begin{equation}
    \label{defGbracketq}
    G_{[q]}(z;z_0,\tau)
    =\sum_{\substack{\ell\le z/\sqrt{q},\\ (\ell,q\tau P(z_0))=1}}
    \frac{\mu^2(\ell)}{\varphi(\ell)}
    \xi_q(z/\ell),
  \end{equation}
  then
  \begin{equation*}
    \xi_q(y)=\sum_{\substack{q_1q_2q_3=q,\\ q_1q_3\le y, \\q_2q_3\le y}}
    \frac{\mu(q_3)\varphi_2(q_3)}{\varphi(q_3)}
    \quad\text{and}
    \quad\varphi_2(q_3)=\prod_{p|q}(p-2).
  \end{equation*}
\end{lem}
\begin{proof}
  We develop the square above and get
  \begin{align*}
    \beta_{z_0,z}(n,\tau)
    &=\sum_{d_1,d_2}\lambda_{d_1}\lambda_{d_2}\1_{[d_1,d_2]|n}
      =\sum_{d_1,d_2}\frac{\lambda_{d_1}\lambda_{d_2}}{[d_1,d_2]}
      \sum_{q|[d_1,d_2}\sum_{a\mode q}e(na/q)
    \\&=\sum_{\substack{q\le z^2,\\ (q,P(z_0))=1}}w_q(z;z_0)c_q(n)
  \end{align*}
  where
  \begin{equation}
    w_q(z;z_0,\tau)
    =
    \sum_{q|[d_1,d_2]}\frac{\lambda_{d_1}\lambda_{d_2}}{[d_1,d_2]}.
  \end{equation}
  We introduce the definition the $\lambda_d$'s, see~\eqref{defbetan},
  and obtain
  \begin{equation*}
    G_\tau(z;z_0)^2
    w_q(z;z_0,\tau)
    =
    \sum_{\substack{\ell_1,\ell_2\le z,\\ (\ell_1\ell_2,\tau P(z_0))=1}}
    \frac{\mu^2(\ell_1)}{\varphi(\ell_1)}
    \frac{\mu^2(\ell_2)}{\varphi(\ell_2)}
    \sum_{\substack{q|[d_1,d_2],\\ d_1|\ell_1, d_2|\ell_2}}\frac{d_1\mu(d_1)d_2\mu(d_2)}{[d_1,d_2]}.
  \end{equation*}
  The inner sum vanishes is $\ell_1$ has a prime factor prime to
  $q\ell_2$, and similarly for $\ell_2$. Furthermore, we need to have
  $q|[\ell_1,\ell_2]$ for the inner sum not be empty. Whence we may
  write $\ell_1=q_1q_3\ell$ and $\ell_2=q_2q_3\ell$ where $(\ell,q)=1$
  and $q=q_1q_2q_3$. The part of the inner sum corresponding to $\ell$
  has value $\prod_{p|\ell}(p-2+1)=\varphi(\ell)$. We have reached
  \begin{equation*}
    G_\tau(z;z_0)^2
    w_q(z;z_0,\tau)
    =
    \mkern-22mu
    \sum_{\substack{\ell\le z,\\ (\ell,q\tau P(z_0))=1}}
    \mkern-17mu
    \frac{\mu^2(\ell)}{\varphi(\ell)}
    \sum_{\substack{q_1q_2q_3=q,\\ q_1q_3\ell\le z, \\q_2q_3\ell\le z}}
    \frac{1}{\varphi(q)\varphi(q_3)}
    \sum_{\substack{q|[d_1,d_2],\\ d_1|q_1q_3,\\ d_2|q_2q_3}}
    \mkern-15mu
    \frac{d_1\mu(d_1)d_2\mu(d_2)}{[d_1,d_2]}.
  \end{equation*}
  In this last inner sum, we have necessarily $d_1=q_1d'_1$ and
  $d_2=q_2d'_2$, so $q_3=[d'_1,d'_2]$. Here is the expression we have
  obtained
  \begin{equation*}
    G_\tau(z;z_0)^2
    w_q(z;z_0,\tau)
    =
    \mkern-22mu
    \sum_{\substack{\ell\le z,\\ (\ell,q\tau P(z_0))=1}}
    \mkern-17mu
    \frac{\mu^2(\ell)}{\varphi(\ell)}
    \sum_{\substack{q_1q_2q_3=q,\\ q_1q_3\ell\le z, \\q_2q_3\ell\le z}}
    \frac{\mu(q)\mu(q_3)}{\varphi(q)\varphi(q_3)}
    \sum_{\substack{q_3=[d'_1,d'_2]}}
    \mkern-15mu
    \frac{d'_1\mu(d'_1)d'_2\mu(d'_2)}{[d'_1,d'_2]}.
  \end{equation*}
  This last inner sum has value $\varphi_2(q_3)$, whence
  \begin{equation*}
    G_\tau(z;z_0)^2
    w_q(z;z_0,\tau)
    =
    \frac{\mu(q)}{\varphi(q)}
    \sum_{\substack{\ell\le z,\\ (\ell,q\tau P(z_0))=1}}
    \frac{\mu^2(\ell)}{\varphi(\ell)}
    \sum_{\substack{q_1q_2q_3=q,\\ q_1q_3\ell\le z, \\q_2q_3\ell\le z}}
    \frac{\mu(q_3)\varphi_2(q_3)}{\varphi(q_3)}
  \end{equation*}
  as announced. The size conditions are readily seen to imply that
  $\ell\le z/\sqrt{q}$.
\end{proof}
\paragraph*{Remark I}
We have
\begin{equation}
  \label{eq:20}
  \left\{
    \begin{aligned}
      &\xi_q(y)=\frac{q}{\varphi(q)}\quad \text{when $y\ge q$},\\
      &|\xi_q(y)|\le \prod_{p|q}\frac{3p-4}{p-1}\quad \text{for every $y>0$}.
    \end{aligned}
    \right.
  \end{equation}
\paragraph*{Remark II}
When developing the theory of \emph{local
  models} in \cite[Chapters 8-11]{Ramare*06}, we show the much
neater expression (see Eq. (11.30) in \cite{Ramare*06}):
\begin{equation*}
  \alpha(n;\tau)=\sum_{d|n}\lambda_d
  =
  \frac{1}{G_\tau(z;z_0)}\sum_{\substack{q\le z,\\
      (q,\tau P(z_0))=1}}\frac{\mu(q)}{\varphi(q)}
  c_q(n).
\end{equation*}
In the theory of local models, we realize $\alpha(n)$ as being (close
to) the best \emph{approximation} of the characteristic function of the
primes, while in the theory of the Selberg sieve, we also ask for a
pointwise \emph{upper bound}.
\paragraph*{Remark III}
In \cite{Hardy*21}, G.H.~Hardy proved in
particular that
\begin{equation}
  \label{eq:2}
  \forall n>1,\quad
  \Lambda(n)=\frac{n}{\varphi(n)}\sum_{\substack{q\ge1}}\frac{\mu(q)}{\varphi(q)}
  c_q(n).
\end{equation}
A proof of this result may be found in the paper \cite{Murty*13} by R.~Murty.
We stress that this expression is \emph{not} valid at
$n=1$. Similarly, if we were to consider the expression
\begin{equation*}
  \sum_{\substack{q\ge1\\ (q,P(z_0))=1}}\frac{\mu(q)}{\varphi(q)}
  c_q(n),
\end{equation*}
we should restrict our attention to integers $n$ that have at least a
prime factor larger than~$z_0$.
\paragraph*{Remark IV}

As can be guessed from
Lemma~\ref{Fourierbetan} and up to some
renormalisation, the coefficients $w_d$ are simply the ``Fourier'' coefficients of the (weighted) sequence
$\beta_{z_0,z}(n)$, as shown in the next easily established expression
(valid for any~$a$ coprime with~$d$):
  \begin{equation*}
    G(z;z_0,\tau)w_d(z;z_0,\tau)=\lim_{N\goes \infty}
    \frac{G(z;z_0,\tau)}{N}\sum_{n\le N}\beta_{z_0,z}(n,\tau)e(na/d).
  \end{equation*}

\section{On the $G$-functions}
\label{Gfunctions}

In this section, we investigate explicit estimates for $G_\tau(z;z_0)$.
The $G$-functions have been thoroughly studied when $z_0=2$ (this
means, no condition on small primes is asked), for instance in
\cite[Lemmas 3.4]{Ramare*95} and more precisely when $\tau=1$ in
\cite{Ramare*14-1}, in \cite{Buethe*14} by Jan B{\"u}the and in
\cite[Theorem 3.1]{Ramare*18-9}.

\begin{lem}
  \label{vanLR}
  When $q$ is prime to $\tau P(z_0)$, we have
  \begin{equation*}
    G_\tau(z/q;z_0)\le \frac{q}{\varphi(q)}G_{q\tau}(z/q;z_0)\le
    G_\tau(z;z_0).
  \end{equation*}
\end{lem}
\noindent
This comes from \cite[Eq. (1.3)]{van-Lint-Richert*65} by J.~van Lint and
H.E.~Richert.
Let us recall \cite[Lemma 2.6]{Ramare*22}.
\begin{lem}\label{EstGdown}
  When $2\le z_0\le z$, we have
  $\displaystyle
  G(z;z_0)\ge e^{-\gamma}\frac{\log z}{\log 2z_0}$.
\end{lem}
We also need an upper estimate in this work.
\cite[Theorem 1.1]{Ramare*14-1} with $q=P(z_0)$ gives a precise answer.
\begin{lem}
  \label{UpperG}
  When $P(z_0)\le z$, $\tau\le z^4$, $(\tau,P(z_0))=1$ and $z_0\ge 35$, we have
  \begin{equation*}
    G_\tau(z^2;z_0)\le 
    3.1\prod_{p< z_0}\frac{p-1}{p}\frac{\varphi(\tau)}{\tau}\log z
    \le 2\frac{\varphi(\tau)}{\tau}\frac{\log z}{\log z_0}.
  \end{equation*}
\end{lem}

\begin{proof}
  \cite[Theorem 1.1]{Ramare*14-1} with $q=\tau P(z_0)$ gives us the bound
  \begin{equation*}
    G_q(z^2)=\frac{\varphi(q)}{q}\biggl(2\log z + \sum_{p|q}\frac{\log
      p}{p-1}+c_0\biggr)
    +\Ocal^*(j(q)/z)
  \end{equation*}
  where
  \begin{equation}
    \label{eq:8}
    j(q)=\prod_{\substack{p|q\\ p\neq
        2}}\frac{p^{3/2}+p-\sqrt{p}-1}{p^{3/2}-\sqrt{p}+1}
    \prod_{\substack{p|q\\ p=2}}\frac{21}{25}
  \end{equation}
  and $c_0$ is given in Lemma~\ref{approxGfunctions}.
  Let us note that
  \begin{equation*}
    \forall p\ge 37,
    \quad
    \frac{p^{3/2}+p-\sqrt{p}-1}{p^{3/2}-\sqrt{p}+1}
    \frac{p}{p-1}\frac{1}{p^{1/20}}\le 1.
  \end{equation*}
  This implies, after a short computation, that $j(\tau P(z_0))\le 15
  (\tau P(z_0))^{1/20}$. We next notice that
  \begin{equation*}
    \sum_{p|q}\frac{\log p}{p-1}
    \le \sum_{p|q}\frac{\log p}{p}+c_0-\gamma\le \log z_0+c_0-\gamma
  \end{equation*}
  by a classical estimate of Rosser \& Schoenfeld in
  \cite[(3.24)]{Rosser-Schoenfeld*62}.
  We have thus reached
  \begin{align*}
    G_\tau(z^2;z_0) = G_q(z^2)
    &\le \prod_{p< z_0}\frac{p-1}{p}\frac{\varphi(\tau)}{\tau}
    \biggl(2\log z +\log z + 2c_0-\gamma + 15\tau^{1/20}/z^{19/20}\biggr)
    \\&\le \prod_{p< z_0}\frac{p-1}{p}\frac{\varphi(\tau)}{\tau}
    \biggl(2\log z +\log z + 2c_0-\gamma + 15/z^{3/4}\biggr)
    \\&\le 3.1\prod_{p< z_0}\frac{p-1}{p}\frac{\varphi(\tau)}{\tau}\log z
  \end{align*}
  since $z\ge P(35)\ge 10^{11}$.
  To handle the contribution of the Euler product over up to~$z_0$, we
  rely on \cite[(3.26)]{Rosser-Schoenfeld*62} and thus
  \begin{equation*}
    \prod_{p< z_0}\frac{p-1}{p}
    \le \frac{e^{-\gamma}}{\log
      z_0}\biggl(1+\frac{1}{2\log^2z_0}\biggr)
    \le \frac{0.585}{\log z_0}.
  \end{equation*}
  The proof is complete.
\end{proof}

\begin{lem}
  \label{Gzz0}
  When $P(z_0)\le z$, $\tau\le z^4$, $(\tau,P(z_0))=1$ and $z_0\ge
  35$, we have
  $\frac{\tau}{\varphi(\tau)}G_\tau(z;z_0)\ge V(z_0)\log(zz_0)
  \ge 0.530\frac{\log z}{\log z_0}$. 
\end{lem}

\begin{proof}
  We use the material of Lemma~\ref{UpperG} to infer that
  \begin{equation*}
    G_\tau(z;z_0)=
    V(z_0)\frac{\varphi(\tau)}{\tau}
    \biggl(\log z + \sum_{p<z_0}\frac{\log p}{p-1}
    +\sum_{p|\tau}\frac{\log p}{p-1}
    +c_0\biggr)
    +\Ocal(15(\tau z)^{1/20}/\sqrt{z}).
  \end{equation*}
  Therefore, and using Lemma~\ref{logz0}, we get
  \begin{equation*}
    G_\tau(z;z_0)
    \ge
    V(z_0)\frac{\varphi(\tau)}{\tau}\biggl(\log (zz_0) -0.6 +c_0-\frac{15}{V(z_0)z^{1/4}}\biggr).
  \end{equation*}
  By Lemma~\ref{getVz0}, we may degrade this inequality in
  \begin{equation*}
    \frac{\tau}{\varphi(\tau)}G_\tau(z;z_0)/V(z_0)\ge
    \log (zz_0) +c_0-0.6-\frac{15 e^\gamma\log(1.13z_0)}{z^{1/4}}.
  \end{equation*}
  Notice that Lemma~\ref{z0vsz} ensures that $\log z\ge \frac45 z_0$,
  so that we only need to check that
  \begin{equation*}
    c_0-0.6-\frac{15 e^\gamma\log(1.13z_0)}{e^{z_0/5}}\ge 0.
  \end{equation*}
  This is readily done (this quantity is even $\ge 0.7$), therefore closing
  the proof of the first inequality. On using Lemma~\ref{getVz0}, we
  further reach the inequality
  \begin{align*}
    \frac{\tau}{\varphi(\tau)}G_\tau(z;z_0)
    &\ge e^{-\gamma}\frac{\log(zz_0)}{\log(1.23 z_0)}
    \ge e^{-\gamma}\frac{\log z}{\log z_0}
    \biggl(\frac{\log z_0}{\log(1.23 z_0)}
    + \frac{\log^2z_0}{\log(1.23 z_0)\log z}\biggr)
      \\&\ge 0.530 \frac{\log z}{\log z_0}
  \end{align*}
  as required.
\end{proof}

\begin{lem}
  \label{Gz2toGz}
  When $P(z_0)\le z$, $\tau\le z^4$, $(\tau,P(z_0))=1$ and $z_0\ge
  35$, we have
  $G_\tau(z^2;z_0)/G_\tau(z;z_0)\le 3.8$. 
\end{lem}

\begin{proof}
  We use Lemmas~\ref{UpperG} and~\ref{Gzz0} to infer that
  \begin{equation*}
    \frac{G_\tau(z^2;z_0)}{G_\tau(z;z_0)}\le \frac{2}{0.530}
    \le 3.78.
  \end{equation*}
  The lemma is proved.
\end{proof}

\begin{lem}
  \label{GztautoGz}
  When $P(z_0)\le z$, $\tau\le z^4$, $(\tau,P(z_0))=1$ and $z_0\ge
  35$, we have
  $G(z;z_0)/G_\tau(z;z_0)\le 1+20/z_0$. 
\end{lem}

\begin{proof}
  By Lemma~\ref{vanLR} (with $q=\tau$ and $\tau=1$),
  we have $G(z;z_0)/G_\tau(z;z_0)\le \varphi(\tau)/\tau$
\end{proof}

\section{On the $w_q$-functions}
\label{wqfunctions}

\subsection{Pointwise estimates}

We follow \cite[Lemma 4.4]{Ramare*95}.
\begin{lem}
  \label{wdforprime}
  When $p$ is prime, we have $-1\le w_p(z;z_0,\tau)G_\tau(z;z_0)(p-1)\le 0$.
\end{lem}

\begin{proof}
  When $p$ is prime and prime to $\tau$, we have $\xi_p(y)=0$ when $y<p$ and
  $\xi_p(y)=p/(p-1)$ otherwise. Hence
  $G_{[p]}(z;z_0)=\frac{p}{p-1}G_{p\tau}(z/p;z_0)$ and the latter quantity
  is at most $G_\tau(z;z_0)$ by Lemma~\ref{vanLR}. The lemma follows readily.
\end{proof}

\begin{lem}
  \label{wdfortwoprimes}
  When $p_1$ and $p_2$ are distinct primes, we have
  \begin{equation*}
    \max\biggl(2,\frac{p_1p_2}{\varphi(p_1p_2)}\biggr)
    \ge w_{p_1p_2}(z;z_0,\tau)G_\tau(z;z_0)\varphi(p_1p_2)
    \ge 0.
  \end{equation*}
\end{lem}

\begin{proof}
  The lemma is trivial if $(p_1p_2,\tau)\neq1$, so let us now assume
  that $p_1$ and $p_2$ are prime to $\tau$.
  Let us assume that $p_1<p_2$. In \eqref{defGbracketq}, when $\ell >
  z/p_2$, there is no solution in $(q_1,q_2,q_3)$. When $z/p_2\ge
  \ell>z/(p_1p_2)$, the only solutions are $(q_1,q_2,q_3)=(p_1,p_2,1)$
  and $(q_1,q_2,q_3)=(p_2,p_1,1)$ so that $\xi_{p_1p_2}(z/\ell)=2$ in
  that range. When $\ell\le z/(p_1p_2)$, we have
  $\xi_{p_1p_2}(z/\ell)=\frac{p_1p_2}{\varphi(p_1p_2)}$. 
\end{proof}

\begin{lem}
  \label{CourageousBis0}
  When $z_0\ge24$, we have
  $\displaystyle
    \bigl|G_\tau(z;z_0)w_q(z;z_0,\tau)\bigr|
    \le 1/q^{2/3}$.

    When $z_0\ge35$, we have
  $\displaystyle
    \bigl|G_\tau(z;z_0)w_q(z;z_0,\tau)\bigr|
    \le 1.04/q^{7/10}$.
\end{lem}

\begin{proof}
  We start from
  \begin{equation*}
    |G_\tau(z;z_0)\varphi(q)w_q(z;z_0,\tau)|
    \le
    \prod_{p|q}\frac{3p-4}{p-1}\frac{G_{q\tau}(z/\sqrt{q}; z_0)}{G_\tau(z;z_0)}
    \le \prod_{p|q}\frac{3p-4}{p-1}.
  \end{equation*}
  We then check that $(3p-4)p^{2/3}/(p-1)^2\le 1$ when $p>23$, hence
  the first result. The second one is proved in much the same manner,
  except that the inequality $f(p)=(3p-4)p^{7/10}/(p-1)^2\le 1$ holds true
  only when $p\ge 43$. We thus have to multiply $1/q^{7/10}$ by
  $f(37)\cdot f(41)$. This proves the second inequality.
\end{proof}

\begin{lem}
  \label{CourageousTer}
  When $z\ge 34$, we have
  \begin{equation*}    
    z_0\max_{q \ge z_0}\bigl|G_\tau(z;z_0)w_q(z;z_0,\tau)\bigr|
    \le
      1+2.2/z_0.
  \end{equation*}
\end{lem}

\begin{proof}
  Any $q$ prime to $P(z_0)$ in the interval $[z_0,z_0^2)$ is a
  prime. By Lemma~\ref{wdforprime}, we then have
  $|G_\tau(z;z_0)w_q(z;z_0,\tau)|\le 1/(q-1)\le 1/(z_0-1)$.
  When $q$ has two prime factors, we use Lemma~\ref{wdfortwoprimes}
  and get
  \begin{equation*}
    G_\tau(z;z_0)w_q(z;z_0,\tau)\le
    \frac{1}{q}\frac{2q}{\varphi(q)}
    \le \frac{2.11}{q}\le
    \frac{2.11}{z_0^2}
    \quad \textsl{(provided that $\omega(q)=2$)}.
  \end{equation*}
  Otherwise, $q\ge
  z_0^3$ and 
  Lemma~\ref{CourageousBis0} tells us that
  $|z_0G_\tau(z;z_0)w_q(z;z_0,\tau)|\le1$. 
  The lemma follows readily from these three observations.
\end{proof}

We may also imitate the proof of \cite[Lemma 4.5]{Ramare*95} to infer
the next lemma.
\begin{lem}
  \label{wdbound}
  We have
   \begin{equation*}
     |w_q(z;z_0,\tau)G_\tau(z;z_0)|\le
     \frac{G_\tau(z^2;z_0)}{G_\tau(z;z_0)}\sum_{\substack{q_1q_2q_3=q\\ 
         q_1q_3,q_2q_3\le z}}\frac{1}{q_1q_2q_3}.
   \end{equation*}
\end{lem}

\begin{proof}
  We easily infer from Lemma~\ref{Fourierbetan} that, when $q$ is
  prime to $\tau P(z_0)$, we have
  \begin{equation*}
    |w_q(z;z_0,\tau)G_\tau(z;z_0)|\le
     \frac{G_{q\tau}(z/\sqrt{q};z_0)}{\varphi(q)G_\tau(z;z_0)}\sum_{\substack{q_1q_2q_3=q\\ 
         q_1q_3,q_2q_3\le z}}1.
   \end{equation*}
   We then complete this inequality with
   $G_{q\tau}(z/\sqrt{q};z_0)/\varphi(q)\le G_\tau(z\sqrt{q};z_0)/q$. 
\end{proof}


  

\subsection{Some variants of the large sieve inequality}

\begin{lem}
  \label{LSspe}
  When $1\le Q_1\le Q_2$ and $\Delta\ge1$ is a positive integer, we
  have
  \begin{equation*}
    \sum_{\substack{Q_1\le q\le Q_2}}\frac{1}{q}
    \sum_{a\mode q\Delta}\biggl|\sum_{n\le N}u_ne\biggl(\frac{na}{q\Delta}\biggr)\biggr|^2
    \le
    (NQ_1^{-1}+2\Delta Q_2)\sum_n|u_n|^2
  \end{equation*}
  for arbitrary complex coefficients $(u_n)$.
\end{lem}
\noindent
\cite[Lemma 8.2]{Ramare*95} already uses a similar mechanism.

\begin{proof}
  Let us first notice that the set of points $(\frac{a}{\Delta q})$
  for $a\mode q\Delta$ and $q\le t$ is $\Delta t^2$-well spaced.
  Let us use the shortcut
  \begin{equation}
    \label{defWq}
    W(q)=\sum_{a\mode q}\biggl|\sum_{n\le N}u_ne(na/q)\biggr|^2
  \end{equation}
    and assume, by homogeneity, that $\sum_n|u_n|^2=1$.
  By
  partial summation, we thus infer that
  \begin{align*}
    \sum_{\substack{Q_1\le q\le Q_2}}\frac{W(\Delta q)}{q}
    &\le \int_{Q_1}^{Q_2}\sum_{q\le t}W(\Delta q)\frac{dt}{t^2}
    +\frac{\sum_{\substack{Q_1\le q\le Q_2}}W(\Delta q)}{Q_2}
    \\&\le \int_{Q_1}^{Q_2}(N+\Delta t^2)\frac{dt}{t^2}
    +\frac{N+\Delta Q_2^2}{Q_2}
    \le \frac{N}{Q_1}+2Q_2,
  \end{align*}
  by using the large sieve inequality on the second line. This ends
  the proof.
\end{proof}

\subsection{A $w_q$-weighted large sieve inequality}
This subsection is devoted to the proof of the next result.
\begin{thm}
  \label{wqLS}
  When $P(z_0)\le z$, $\tau\le z^4$, $(\tau,P(z_0))=1$, $z_0\ge
  35$ and $z_1\in[z_0,\sqrt{z}/\log z]$, we have
  \begin{multline*}
    \sum_{z_1\le  q\le z^2}
    |G_\tau(z;z_0)w_q(z;z_0,\tau)|\sum_{a\mode q}\biggl|\sum_{n\le N}u_ne(na/q)\biggr|^2
    \\\le 13 (Nz_1^{-1}+z^2(\log 3z_0)^{-1}) \sum_{n}|u_n|^2
  \end{multline*}
  for arbitrary complex coefficients $(u_n)$.
\end{thm}
\noindent
This is the extension of \cite[Theorem 8.1]{Ramare*95}, though we did
not try here to optimize the constants; we restricted our attention to
eliminating the log-powers that intervene with a more casual
treatment. 

\begin{proof}[Proof of Theorem~\ref{wqLS}]
  Let us use the notation given by \eqref{defWq}
  and assume, by homogeneity, that $\sum_n|u_n|^2=1$.
  We readily check, via Lemma~\ref{z0vsz}, that the condition $P(z_0)\le z$ implies that
  $z_0\le z^{1/2}/\log z$. 
  When $q\le z$, we may rely on Lemma~\ref{CourageousTer} and write
  (recall $W(q)$ is defined in~\eqref{defWq})
   \begin{align}
    \Sigma(z_0,z)=\sum_{z_0\le  q\le z}
    |G_\tau(z;z_0)w_q(z;z_0,\tau)|W(q)
    &\le\frac{1+2.2/z_0}{z_0}\sum_{z_0\le  q\le z}
      W(q)\notag
     \\&\le \frac{1+2.2/z_0}{z_0}(N+z^2)
     \label{Ineq1}
  \end{align}
  by the large sieve inequality.
  When $q\le z^{3/2}$, we may use Lemma~\ref{CourageousBis0} and write
  \begin{align*}
    \Sigma(z,z^{3/2})&=\sum_{z\le  q\le z^{3/2}}
    |G_\tau(z;z_0)w_q(z;z_0,\tau)|W(q)
    \le 1.04\sum_{z\le  q\le z^{3/2}}
      W(q)/q^{7/10}
    \\&\le 1.04\int_{z}^{z^{3/2}}\sum_{q\le
    t}W(q)\frac{(7/10)dt}{t^{17/10}}
    +1.04\frac{\sum_{q\le z^{3/2}}W(q)}{z^{7/5}}
    \\&\le 1.04 \int_{z}^{z^{3/2}}(N+t^2)\frac{(7/10)dt}{t^{17/10}}
    +1.04\frac{N+z^3}{z^{7/5}}    
  \end{align*}
  by the large sieve inequality. This simplifies into
  \begin{equation}
    \label{Ineq2}
    \Sigma(z,z^{3/2})/1.04
    \le
    \frac{N}{z^{7/10}}+\biggl(\frac{7}{10}\frac{10}{13}+1\biggr)z^{\frac{3}{2}\frac{13}{10}}
    =\frac{N}{z^{7/10}}+\frac{20}{13}z^{2-\frac1{20}}.
  \end{equation}
  Let us now handle the more difficult part. We employ
  Lemma~\ref{wdbound} together with Lemma~\ref{Gz2toGz} to write:
  \begin{multline*}
    \Sigma(z^{3/2},z^2)=\sum_{z^{3/2}\le  q\le z^{2}}
    |G_\tau(z;z_0)w_q(z;z_0,\tau)|W(q)
    \le
      3.8 \sum_{\substack{z^{3/2}\le  q_1q_2q_3\le z^{2}\\
    q_1q_3,q_2q_3\le z\\ (q_1q_2q_3,P(z_0))=1}}
    \frac{W(q_1q_2q_3)}{q_1q_2q_3}
    \\\le 3.8
    \sum_{\substack{q_3\le z\\
    (q_3,P(z_0))=1}}\frac{\mu^2(q_3)}{q_3}
    \sum_{\substack{q_1\le z/q_3\\
    (q_1,q_3P(z_0))=1}}\frac{\mu^2(q_1)}{q_1}
    \sum_{\substack{z^{3/2}/(q_1q_3)\le  q_2\le z/q_3\\
     (q_2,q_1q_3P(z_0))=1}}
    \frac{\mu^2(q_2)W(q_1q_2q_3)}{q_2}.
  \end{multline*}
  Notice that, for the last summation to be non-empty, we should
  require that $q_1\ge z^{1/2}$ and thus $q_3\le z^{1/2}$. We get
  \begin{align*}
    \frac{\Sigma(z^{3/2},z^2)}{3.8}
    &\le
    \sum_{\substack{q_3\le z^{1/2}\\
    (q_3,P(z_0))=1}}\frac{\mu^2(q_3)}{q_3}
    \sum_{\substack{z^{1/2}\le q_1\le z/q_3\\
        (q_1,q_3P(z_0))=1}}\frac{\mu^2(q_1)}{q_1}
    \biggl(\frac{q_1q_3N}{z^{3/2}}+2q_1q_3(z/q_3)\biggr)
    \\&\le
    \frac{N}{z^{3/2}}\sum_{\substack{q_3\le z^{1/2}\\
    (q_3,P(z_0))=1}}
    \sum_{\substack{ q_1\le z/q_3\\ (q_1,P(z_0))=1}}1
    + 2z \sum_{\substack{q_3\le z^{1/2}\\
    (q_3,P(z_0))=1}}\frac{\mu^2(q_3)}{q_3}
    \sum_{\substack{ q_1\le z/q_3\\ (q_1,P(z_0))=1}}1
    \\&\le
    1.1\frac{N}{z^{1/2}}\frac{\log z}{\log 3z_0}
    + \frac{3z^2}{\log (3z_0)}
  \end{align*}
  by using Lemma~\ref{LSspe} on the first line and
  Lemma~\ref{bettersum} together with the easily
  established inequalities
  \begin{equation*}
    \sum_{\substack{q_3\le z^{1/2}\\
    (q_3,P(z_0))=1}}\frac{\mu^2(q_3)}{q_3}\le \log z,\quad
    \prod_{p\ge z_0}\biggl(1+\frac{1}{p^2}\biggr)\le \frac{3}{2\cdot 1.1}.
  \end{equation*}
  on the third line.
  On collecting our estimates, we finally find that
  \begin{align*}
    \Sigma(z_1,z^2)
    &\le \Sigma(z_1,z)+\Sigma(z,z^{3/2})+\Sigma(z^{3/2},z^2)
    \\&\le \frac{N}{z_1}\biggl(1+\frac{2.2}{35}
    +1.04\frac{z_1}{z^{7/10}}+3.8\times1.1\frac{z_1\log z}{z^{1/2}\log(3z_0)}\biggr)
    \\&\qquad+\frac{z^2}{\log
    3z_0}\biggl(\frac{1+\frac{2.2}{35}}{z_0}\log(3z_0)
    +1.04\frac{7\log(3z_0)}{13 z^{1/20}}+3.8\times3\biggr).
  \end{align*}
  We use $z\ge \exp(4z_0/5)$ given by Lemma~\ref{z0vsz} to show that
  \begin{equation*}
    \left\{\begin{aligned}
        &1+\frac{2.2}{35}
    +1.04\frac{z^{1/2}}{z^{7/10}\log z}+\frac{3.8\times1.1}{\log(3z_0)}
    \le 1.97\\
    &\frac{1+\frac{2.2}{35}}{z_0}\log(3z_0)
    +1.04\frac{20\log(3z_0)}{13 z^{1/20}}+3.8\times3
    \le 12.2.
  \end{aligned}
      \right.
  \end{equation*}
  The theorem is proved.
\end{proof}

\section{A large sieve inequality with primes and few points}

\begin{lem}
  \label{BeurlingF}
  Let $M\in\mathbb{R}$, and $N$ and $\delta$ be positive real number.
  There exists a smooth function $\psi$ on $\mathbb{R}$ such that
  \begin{itemize}
  \item The function $\psi$ is non-negative.
  \item When $t\in [M,M+N]$, we have $\psi(t)\ge 1$.
  \item $\psi(0)=N+\delta^{-1}$.
  \item When $|\alpha|>\delta$, we have $\hat{\psi}(\alpha)=0$.
  \item We have $\psi(t)=\Ocal_{M,N,\delta}(1/(1+|t|^2))$.
  \end{itemize}
\end{lem}
This is \cite[Lemma 3.3]{Ramare*22} and is a retelling of a result due
to A. Selberg, see \cite[Section 20]{Selberg*91}.

\begin{thm}
  \label{Bel} Let $N\ge 10^4$.
  Let $\mathcal{B}$ be a $\delta$-well spaced subset of
  $\mathbb{R}/\mathbb{Z}$. For any function $f$
  on $\mathcal{B}$, we have
  \begin{equation*}
    \sum_{p\le N}\biggl|\sum_{b\in \mathcal{B}}f(b)e(bp)\biggr|^2
    \le 19
    (N+\delta^{-1})\|f\|_2^2\frac{\log
      (2\|f\|_1^2/\|f\|_2^2)}{\log N}.
  \end{equation*}
  where $\|f\|_q^q=\sum_{b\in \mathcal{B}}|f(b)|^q$ for any positive $q$.
\end{thm}
\noindent
The best constant that comes out of the proof we propose is
$(4e^\gamma+o(1))$, provided~$\|f\|_1^2/\|f\|_2^2$ goes to infinity.
\begin{proof}
  Let us first notice that $\|f\|_1^2\ge\|f\|_2^2$.
  \subsubsection*{Small $\|f\|_1^2/\|f\|_2^2$}
  When $y=\|f\|_1^2/\|f\|_2^2$ is small, we simply use
  \begin{equation*}
    \sum_{p\le N}\biggl|\sum_{b\in \mathcal{B}}f(b)e(bp)\biggr|^2
    \le \pi(N)\|f\|_1^2\le \frac{5N}{4\log N}\|f\|_2^2
    \frac{y}{\log(2y)}\log
    (2\|f\|_1^2/\|f\|_2^2)
  \end{equation*}
  by Lemma~\ref{RS} and provided that $N\ge 114$. When $y\le 21$, we obtain
  \begin{equation*}
    \sum_{p\le N}\biggl|\sum_{b\in \mathcal{B}}f(b)e(bp)\biggr|^2
    \le 4e^\gamma \|f\|_2^2\frac{\log
      (2\|f\|_1^2/\|f\|_2^2)}{\log N}.
  \end{equation*}
  \subsubsection*{Large $\|f\|_1^2/\|f\|_2^2$}
  When $\|f\|_1^2/\|f\|_2^2\ge N^{e^{-\gamma}/4}$, we use the dual of the usual large sieve
  inequality (see \cite{Montgomery*78} by H.L. Montgomery) to infer that
  \begin{align*}
    \sum_{p\le N}\biggl|\sum_{b\in \mathcal{B}}f(b)e(bp)\biggr|^2
    &\le
      (N+\delta^{-1})\|f\|_2^2
    \le  (N+\delta^{-1})\|f\|_2^2\frac{\log
    (2\|f\|_1^2/\|f\|_2^2)}{\log (N^{e^{-\gamma}/4})}
    \\&\le
    4e^\gamma (N+\delta^{-1})\|f\|_2^2\frac{\log
    (2\|f\|_1^2/\|f\|_2^2)}{\log N}.
  \end{align*}
  This establishes our inequality
  in this case.
  \subsubsection*{Small primes}
  Let $z=N^{1/4}/10$ and
  \begin{equation*}
    z_0=\frac{2\|f\|_1^2}{1.23\|f\|_2^2}> 34.
  \end{equation*}
  We assume that $21<\|f\|_1^2/\|f\|_2^2\le N^{e^{-\gamma}/4}$ and
  that $N\ge 10^4$.
  We discard the small primes trivially:
  \begin{align*}
    \sum_{ p\le z}\biggl|\sum_{b\in \mathcal{B}}f(b)e(bp)\biggr|^2
    &\le z\|f\|_1^2\le N^{(1+e^{-\gamma})/4}\|f\|_2^2/10
    \\&\le N \frac{\|f\|_2^2\log(2{\|f\|_1^2}/{\|f\|_2^2})}{\log
    N}
    \frac{\log N}{10N^{(3-e^{-\gamma})/4}\log 42}
    \\&\le \frac{N}{1113} \frac{\|f\|_2^2\log(2{\|f\|_1^2}/{\|f\|_2^2})}{\log
    N}.
  \end{align*}
  \subsubsection*{Main proof}
  Let us define
  \begin{equation}
    \label{defW}
    W=\sum_{z<p\le N}\biggl|\sum_{b\in \mathcal{B}}f(b)e(bp)\biggr|^2.
  \end{equation}
  We bound above the characteristic function of the primes from $(z,N]$ by
  our enveloping sieve and further majorize the characteristic
  function of the interval~$[1,N]$ by a function $\psi$ (see
  Lemma~\ref{BeurlingF}) of Fourier transform supported by~$[-\delta_1,\delta_1]$ where $\delta_1=\min(\delta,1/(2z^4))$, and
  which is such that $\hat{\psi}(0)=N+\delta_1^{-1}$. This leads to
  \begin{equation*}
    W\le
    \sum_{\substack{q\le z^2,\\
        (q,P(z_0))=1}}w_q(z;z_0)
    \sum_{a\mode  q}\sum_{b_1,b_2}f(b_1)\overline{f(b_2)}
    \sum_{n\in\mathbb{Z}}e((b_1-b_2)n)e(an/q)\psi(n).
  \end{equation*}
  We have shorthened $w_q(z;z_0,1)$ in $w_q(z;z_0)$.
  We split this quantity according to whether $q< z_0$ or not:
  \begin{equation*}
    W=W(q< z_0)+ W(q\ge z_0).
  \end{equation*}
  When $q\ge z_0$,
  Poisson summation formula tells us that the inner sum
  reads as $\sum_{m\in\mathbb{Z}}\hat{\psi}(b_1-b_2-(a/q)+m)$.
  The sum over $b_1$, $b_2$ and $n$ is thus
  \begin{equation*}
    \le (N+\delta^{-1}_1)\sum_{b_1,b_2}|f(b_1)||f(b_2)|\#\bigl\{a/q:
     \|b_1-b_2+a/q\|< \delta_1
    \bigr\}.
  \end{equation*}
  Given $(b_1,b_2)$ and since
  $1/z^4>2\delta_1$, at most one $a/q$ may work. On bounding above~$|w_q(z;z_0)|$ by
  Lemma~\ref{CourageousTer}, we see that
  \begin{equation}
    \label{MainStep}
    G(z;z_0)W(q\ge z_0)
    \le
    (1+2.2\,z_0^{-1})(N+\delta^{-1}_1)\frac{\|f\|_1^2}{z_0}.
  \end{equation}
  When $w_q(z;z_0)\neq0$, we have $q|P(z)/P(z_0)$; on adding the
  condition~$q<z_0$, only $q=1$ remains.
  Since $\mathcal{B}$ is $\delta$-well-spaced and $w_1(z;z_0)=1/G(z;z_0)$,
  we infer that
  \begin{equation*}
   G(z;z_0)W(q<z_0)
    \le
    (N+\delta_1^{-1})\|f\|_2^2.
  \end{equation*}
  On recalling Lemma~\ref{EstGdown}, we thus get
  \begin{equation}
    \label{eq:1}
    W \le (N+\delta^{-1}_1)
    \biggl(\|f\|_2^2+(1+2.2\,z_0^{-1})\frac{\|f\|_1^2}{z_0}\biggr)
    \frac{e^\gamma\log (1.23\,z_0)}{\log z}.
  \end{equation}
  \subsubsection*{Final estimate}
  We check that $(N+\delta_1^{-1})\le
  \frac{N+4N10^{-4}}{N}(N+\delta^{-1})$. We finally get
  \begin{multline*}
    \sum_{ p\le N}\biggl|\sum_{b\in \mathcal{B}}f(b)e(bn)\biggr|^2\le
    \biggl(\frac1{1113}
    +
    4e^\gamma(1+4/10^4)\biggl(1+(1+2.2/34)\frac{1.23}{2} \biggr)
    \frac{\log N}{\log N-4\log 4}\biggr)
    \\\times(N+\delta^{-1})\|f\|_2^2\frac{\log
      (2\|f\|_1^2/\|f\|_2^2)}{\log N}.
  \end{multline*}
  The proof of the theorem follows readily.
\end{proof}

\section{Proof of Theorem~\ref{ExtensionpreciseNow}} 

\begin{proof}[Proof of Theorem~\ref{ExtensionpreciseNow}]
  We follow the usual duality argument, starting from Theorem~\ref{Bel}. We
  write
  \begin{equation*}
    \sum_{x\in\Xcal}
    \biggl|\sum_{ p\le N}{}u_p e(xp)\biggr|^2
    =\sum_{p\le N}u_p\sum_{x\in\Xcal}\overline{S(x)}e(xp)
  \end{equation*}
  with $S(x)=\sum_{ p\le N}{}u_p e(xp)$. We apply Cauchy's inequality
  to the resulting expression, then Theorem~\ref{Bel} and the
  inequality
  \begin{equation*}
    \biggl(\sum_{x\in\Xcal}|S(x)|\biggr)^2/\sum_{x\in
      \Xcal}|S(x)|^2\le |\Xcal|
  \end{equation*}
  to simplify the factor $\log(2\|f\|_1^2/\|f\|_2^2)$ that appears in
  Theorem~\ref{Bel}. The inequality of the theorem then follows
  swiftly. On taking $u_p=1$ and $\Xcal=\{a/q:q\le Q,a(q)=1\}$ where
  we have $S(a/q)=\sum_{p\le N}e(ap/q)=\mu(q)(1+o(1))N/(\varphi(q)\log N)$
  when $Q$ is at most a power of $\log N$, we obtain that the quantity
  $\sum_{x\in\Xcal}|S(a/q)|^2$ is asymptotic to $G(Q)\pi(N)^2$. The
  last part of the theorem then follows on noticing that $\log
  |\Xcal|\sim2\log Q\sim 2G(Q)$.
\end{proof}

\part{Cusps}

\section{Bounding the number of cusps}
We define, as in Eq.\,(14) in \cite{Ramare*22},
\begin{equation}
  \label{defV}
  V = \biggl(\frac{N+\delta^{-1}}{\log N}\sum_{p\le
    N}|{}u_p|^2\biggr)^{1/2}.
\end{equation}
\begin{lem}
  \label{SizeLevelSets}
  Let $N\ge 10^4$ and $\Xcal$ be a $\delta$-well spaced subset of
  $\mathbb{R}/\mathbb{Z}$. Let $({}u_p)_{p\le N}$ be a sequence of complex numbers.
  We have, for $A\ge 1$,
  \begin{equation*}
    \#\biggl\{x\in\Xcal:\biggr|\sum_{ p\le N}{}u_p e(xp)\biggl|\ge V/A\biggr\}
    \le
    19A^2\log (2A),
  \end{equation*}
  where $V$ is defined in \eqref{defV}.
  The constant 19 may be replaced by $4e^\gamma+o(1)=7.12\cdots$ when
  $N\goes\infty$ and $A\goes \infty$.
\end{lem}
We used the notation $\#S$ to denote the cardinality of the
set~$S$. 
\begin{proof}[Proof of Lemma~\ref{SizeLevelSets}]
  This is a trivial applications of Markov's inequality.
\end{proof}

\begin{proof}[Proof of Theorem~\ref{L1}]
  We readily
  getting
  \begin{align*}
    \int_0^1 |T^*(\alpha)|d\alpha
    &\ll
    \sum_{0\le r\le R-1}\frac{e^{2r}K}{N}f(r)\frac{T^*(0)}{e^r}
    +\frac{T^*(0)}{e^{R}}
    \\&\ll
    T^*(0)\biggl(\frac{e^{R}K}{N}f(R)
    +e^{-R}\biggr)
    \ll
    \frac{\sqrt{Nf(\log N)}}{\sqrt{K}\log N}
  \end{align*}
  by choosing $R=[(1/2)\log \frac{N}{Kf(\log N)}]$. The proof is complete.
\end{proof}

\begin{proof}[Proof of Theorem~\ref{thmCusps}]
Lemma~\ref{SizeLevelSets} tells us, after a change of
variable ($A\mapsto A\sqrt{K(1+(N\delta)^{-1})}$), that
\begin{equation}
  \label{DensityPrimes}
  \#\biggl\{x\in\Xcal:\biggr|\sum_{\substack{p\le N\\
      p\in\mathcal{P}^*}}
  e(xp)\biggl|\ge \frac{N}{AK\log N} \biggr\}
  \le
  19A^2K(1+(N\delta)^{-1}).
\end{equation}
The above discussion leads to Theorem~\ref{thmCusps}.
\end{proof}

\section{Finding rational points as cusps}
\label{LP}
\subsubsection*{On the structure of the set of cusps}
Lemma~\ref{thmCusps} gives an upper bound for the number of cusps. Let
us now investigate the structure of this set. The examples we
provided tell us it is possible that this set has a structure close to
$\mathcal{F}+\bigl\{\frac{a}{q},q\le 100A\bigr\}$ for a small enough
set~$\mathcal{F}$. This is also a rationale leading to the proof of
Theorem~\ref{ExtensionpreciseNow} and why we thought interesting not
to sieve the small moduli.  Here is a first lemma.
\begin{lem}
  \label{iniStructureCusps}
  Notation being as in Definition~\ref{defCusps}, and assuming
  $N\ge 10$, the following holds:
  \begin{itemize}
  \item When $\alpha$ lies in
    $\Cscr(\mathcal{P}^*,A)$, so do $-\alpha$ and
    $\frac12+\alpha$.
  \item When
    $\{\frac13,\frac12,\frac23,1\}\subset\Cscr(\mathcal{P}^*,1/2)$ 
  \item For every $\xi\in\RsurZ$ with $T^*(\xi)\neq0$ and every
    square-free positive integer $q<\sqrt{N}$ such that
    $\varphi(q)\le AT^*(0)/|T^*(\xi)|$, there exists $a$, coprime with $q$,
    such that $\xi +(a/q)$ lies in $\Cscr(\mathcal{P}^*,A)$.
  \end{itemize}
\end{lem}

\begin{proof}
  Let us recall that $\mathcal{P}^*\in[\sqrt{N},N]$, so that the
  elements of $\mathcal{P}^*$ are prime to the moduli 2, 3 and $q$
  appearing in the three claimed properties.
  The first item is a consequence
  of the facts that (1) the characteristic function
  $\1_{\mathcal{P}^*}$ of~$\mathcal{P}^*$ is real valued and (2) that every
  member of $\mathcal{P}^*$ is odd.
  Concerning the third item,
  it is enough to notice the inequality
  \begin{equation}
    \label{condgen}
    \sum_{a\mode q}
    \biggr|\sum_{\substack{p\in\mathcal{P}^*}}
    e(p(\xi+ (a/q))\biggl|
    \ge \biggr|\sum_{a\mode q}
    \sum_{\substack{p\in\mathcal{P}^*}}
    e(p(\xi+ (a/q))\biggl|=\mu^2(q)|T^*(\xi)|
  \end{equation}
  with $\xi=0$, since $c_q(p)=\mu(q)$ where $c_q(n)$ denotes the
  value of the Ramanujan sum at~$n$.
  The second item follows from this same inequality applied to $q=3$,
  on noticing that the absolute values of the involved Fourier
  polynomial at $1/3$ and $2/3$ are the same.
\end{proof}

The proof of Theorem~\ref{Companions} relies on extracting information from the stream of
inequalities~\eqref{condgen}, for varying $q$'s. 
\begin{proof}[Proof of Theorem~\ref{Companions}]
  We may assume that $B\le A/10$, for otherwise
  the result is obvious: $\mathcal{F}$ contains the element~$\xi$.
  Let us set $Q=A/(2B)\ge1$ and
  \begin{equation}
    \label{eq:9}
    \mathcal{Q}=\bigl\{q: q\le Q,\mu^2(q)=1\bigr\}.
  \end{equation}
  The non-negative real variables
  $t_\xi(a/q)=t^*(\xi+(a/q))$ are
  bounded above by~$Z\le 1$ and  satisfy (see~\eqref{condgen})
  \begin{equation*}
    \forall q\in\mathcal{Q},
    \quad
    \sum_{a\mode q}t_\xi(a/q)\ge 1/B
  \end{equation*}
  as well as, by Lemma~\ref{thmCusps} and since the points in
  $\mathcal{F}$ are $1/N$-well spaced,
  \begin{equation*}
    \forall C\ge 1, \quad
    \#\bigl\{a/q:q\in\mathcal{Q}, t_\xi(a/q)\ge 1/C\bigr\}\le 19
    KC^2\log 2C.
  \end{equation*}
  Let us further introduce the variables
  \begin{equation}
    \label{eq:10}
    x(q,C)=\#\{a\mode q: bZ/C > t_\xi(a/q) \ge bZ/(C+1)\}
  \end{equation}
  for some $b\in(1,2)$.
  We put $bZ/C$, rather than $Z/C$, so as to include the case $t_\xi(a/q)=Z$
  in the same setting. 
  We thus get
  \begin{equation*}
    \forall q\in\mathcal{Q},
    \quad
    \sum_{C\ge 1}\frac{bZ}{C}x(q, C)\ge 1/B.
  \end{equation*}
  Let us assume momentarily that $AZ=A'$ is a positive integer and degrade the above
  inequality into
  \begin{equation*}
    \sum_{1\le C\le A'-1}\frac{b}{C}x(q, C)
    +\Bigl(\varphi(q)-\sum_{1\le C\le A'-1}x(q, C)\Bigr)\frac{b}{A'}
    \ge 1/(BZ),
  \end{equation*}
  i.e. 
  \begin{equation*}
    \sum_{1\le C\le A'-1}\biggl(\frac{b}{C}-\frac{b}{A'}\biggr)x(q, C)
    \ge \frac{1}{BZ}-\frac{b\varphi(q)}{A'}\ge 
    \frac{1-b/2}{BZ}.
  \end{equation*}
  We now sum over $q\in\mathcal{Q}$,
  getting, by Lemma~\ref{SQp},
  \begin{equation*}
    \sum_{1\le C\le A'-1}\biggl(\frac{b}{C}-\frac{b}{A'}\biggr)
    \sum_{q\in\mathcal{Q}}x(q, C)
    \ge \frac{(1-b/2)A'}{4B^2Z^2}.
  \end{equation*}
  Let us set
  \begin{equation*}
    X(C)=\sum_{D\le C}\sum_{q\in\mathcal{Q}}x(q, D),\quad
    (X(0)=0)
  \end{equation*}
  which satisfies
  \begin{equation*}
    \sum_{1\le C\le A'-1}\biggl(\frac{b}{C}-\frac{b}{A'}\biggr)
    (X(C)-X(C-1))
    \ge \frac{(1-b/2)A'}{4B^2Z^2}.
  \end{equation*}
  On reshuffling the left-hand side, we obtain:
  \begin{equation*}
    \sum_{1\le C\le A'-1}X(C)\biggl(\frac{b}{C}-\frac{b}{C+1}\biggr)
    +X(A'-1)\biggl(\frac{b}{A'-1}-\frac{b}{A'}\biggr)
    \ge
    \frac{(1-b/2)A'}{4B^2Z^2}.
  \end{equation*}
  The non-decreasing function $C\mapsto X(C)$ is therefore constrainted by
  the two inequalities:
  \begin{equation*}
    \sum_{1\le C\le A'-1}\frac{X(C)}{C(C+1)} \ge
    \frac{(1-b/2)A'}{4B^2Z^2},
    \quad
    X(C)\le 19 K(C+1)^2\log(2C+2)/(bZ)^2.
  \end{equation*}
  Let us split the first sum at $C=\theta A'$. We deduce from the
  above that
  \begin{equation*}
    19K(bZ)^{-2}\log(\theta A'+1)\sum_{1\le C< \theta A'}\frac{C+1}{C}
    +X(A'-1)\biggl(\frac{1}{[\theta A']}-\frac{1}{A'}\biggr)
    \ge
    \frac{(1-b/2)A'}{4B^2Z^2}
  \end{equation*}
  and thus
  \begin{equation*}
    38K(bZ)^{-2}\theta A'\log(2\theta A'+2)
    +X(A'-1)\frac{4-\theta+A^{\prime-1}}{4\theta A'}
    \ge \frac{(1-b/2)A'}{4B^2Z^2}.
  \end{equation*}
  This calls for $\theta \log A'=b^2/(8\times38\times2 B^2K)$, so that
  \begin{equation*}
    \tfrac98 X(A'-1)\ge (3-2b)\frac{A^{\prime2}b^2}{4864
      K(BZ)^2B^2\log 2A'}.
  \end{equation*}
  On letting $b$ go to 1, this finally amounts to
  $
    X(A'-1)\ge  \frac{A^2}{5472Z^2B^4K\log 2A'}
    $. This is proved when $AZ$ is an integer. Otherwise, we use the
    same bound but for $A"=[AZ]/Z\ge A-1/Z\ge (1-1/(ZA))A\ge 9A/10$    
  to conclude the proof.
\end{proof}

\part{The structure theorem}
\label{Transference}

\newcommand{\myxi}{y}

\section{Preliminaries}
It is better for clarity to present the parameters first and the
object to study later.
So let $\epsilon\in(0,1/2]$ and a parameter
$A\ge1$ be given.
We define the parameter $z_0$ by:
\begin{equation}
  \label{defz0}
  z_0=\exp\bigl(25000A^3K(\log 2A)^2\bigr)/\epsilon.
\end{equation}
(Notice for numerical purposes that $z_0\ge 35$)
The reader may want to keep this parameter as unknown until
Theorem~\ref{MainDec}, so as to understand the above choice.
We further introduce a positive integer parameter $M$ that satisfies
\begin{equation}
  \tag{H1}
 \text{Every prime $p|M$ is $<z_0$ and $P(z_0)|M$.}
\end{equation}
At first, the readers may select $M=P(z_0)$, but an application we
have in mind
will require $M=P(z_0)^3$.

Let then  $\mathcal{P}^*$ be a subset of the primes within $[1,N]$, of
cardinality $N/(K\log N)$, all larger than $z=\sqrt{N/(Mz_0)}$.
We assume that $P(z_0)\le z$. By \cite[Theorem
9]{Rosser-Schoenfeld*62} of J.-B. Rosser and L. Schoenfeld, this is
implied by $1.02 z_0\le \log z$.
\begin{equation}
  \tag{H2}
  P(z_0)\le z\quad\text{and}\quad
  z\ge N^{1/{2(1+\epsilon)}}.
\end{equation}
The last inequality is equivalent to $Mz_0\ge N^{\frac{\ve}{1+\ve}}$.
We define
\begin{equation}
  \label{defNprime}
  T^*(\alpha)=\sum_{p\in\mathcal{P}^*}e(p\alpha),
  \quad
  N'=240 NA,
  \quad
  \varepsilon=\frac{1}{240A}.
\end{equation}
We denote by $f$ the characteristic function of $\mathcal{P}^*$.

\section{A finite cover of $\Cscr(\mathcal{P}^*,A)$}
\label{Cover}
Let us cover the set $\Cscr(\mathcal{P}^*,A)$ of $A$-cusps by a finite
set $\Xi$ of points so that every points of $\Cscr(\mathcal{P}^*,A)$
is at distance $\le 1/(240A)$ of a point of $\Xi$.
To do so, we take, if possible, a point $\myxi$ in each interval
$[\frac{a-1}{N'},\frac{a}{N'})$ with $a\le N'=240AN$ and even (resp. odd)
such that
\begin{equation}
  \label{eq:6}
  |T^*(\myxi)|\ge T^*(0)/A.
\end{equation}
By Eq.~\eqref{DensityPrimes} with $\delta=1/N'$, each set (with $a$ odd and with $a$
even) has at most $5000A^3K\log(2A)$ points. The set $\Xi$ is the union of
both. Every point of $\Cscr(\mathcal{P}^*,A)$ is at a distance
$\le \ve/N$ of a point in $\Xi$.  We also consider
$\Xi_M=\{M\myxi:\myxi\in\Xi\}\subset\RsurZ$, which may be much smaller
that $\Xi$ (see Theorem~\ref{Companions}), it may even be reduced to one point.
\section{The associated Bohr set}
We consider
\begin{equation}
  \label{defBM}
  \Bcal_M(\ve)=\bigl\{n\le N: M|n\quad \&\quad
  \forall \myxi\in\Xi, \|\myxi n\|\le \ve\bigr\},
\end{equation}
as well as
\begin{equation}
  \label{defSM}
  S_M(\alpha; \Xi,\ve)=\sum_{b\in\Bcal_M(\ve)}e(b\alpha).
\end{equation}
See \cite[Proposition 4.2]{Green-Ruzsa*07} by B.~Green \& I.~Ruzsa
or \cite[Lemma 10.4]{Green-Tao*08-2}.

The parameter $N$ is not recalled in our notation.  By \cite[Lemma
4.20]{Tao-Vu*06}\footnote{A more usual argument based on the pigeonhole
  principle yields the marginally weaker bound
  $\ve^{|\Xi|}N/(1+\ve)^{|\Xi|}$. Asymptotically, when $N$ goes to
  infinity, $\ve$ remains fixed and the set $\Xi$ is made of
  $\mathbb{Q}$-linear independant points, we have
  $|\Bcal|\sim (2\ve)^{|\Xi|}N$.}  of the book of T.~Tao \& V.~Vu,
applied to $n/M$ and $\Xi_M$, we have
\begin{equation}
  \label{sizeB}
  |\Bcal_M|\ge \ve^{|\Xi_M|}\lfloor N/M\rfloor
  \ge
  \tfrac12 \ve^{|\Xi_M|}N/M
  .
\end{equation}

\begin{lem}
  \label{tutu}
  When $u\in\RsurZ$, we have $|e(u)-1|\le 2\pi \|u\|$, where $\|u\|$
  is the distance to the nearest integer.
\end{lem}

\begin{proof}
  First replace $u$ by $w=\pm u+k$ for some $k\in\mathbb{Z}$, in
  such a way that $w=\{u\}$. 
  The result then comes from the mean value theorem or from
  $e(w)-1=\int_0^w 2i\pi e(v)dv$. 
\end{proof}

\begin{lem}
  \label{GoodPol}
  If, given $\alpha\in\RsurZ$, there exists $\myxi\in\Xi$ such that $|\alpha-\myxi|\le \ve /N$,
  then $S_M(\alpha; \Xi,\ve) =(1+\Ocal^*(7\ve))|\Pcal_M(\ve)|$.
\end{lem}

\begin{proof}
  For any $b\in\Bcal_M(\ve)$, we have $\|b\alpha\|\le
  2\ve$. Therefore, on calling to Lemma~\ref{tutu}, we get
  \begin{equation*}
    S_M(\alpha; \Xi,\ve)
    =\sum_{b\in\Pcal_M(\ve)}\bigl(1+e(b\alpha)-1\bigr)
    =(1+\Ocal^*(7\ve))|\Bcal_M(\ve)|.
  \end{equation*}
  Once this is noticed, the lemma follows readily.
\end{proof}
Furthermore and following B. Green in \cite[Eq. (6.5)]{Green*05}, we then consider
\begin{equation}
  \label{defrhoM}
  \rho_{M,\ve}(m)=\frac{1}{|\Bcal_{M}(\ve)|^2}
  \sum_{\substack{b_1-b_2=m\\ b_1,b_2\in\Bcal_M(\ve)}}\mkern-20mu1\mkern20mu
  \ge0.
\end{equation}

\section{The approximant}
We consider
\begin{equation}
  \label{deffast}
  \begin{array}{rcl}
    f^\ast:\mathbb{Z}&\rightarrow&\mathbb{C}\\
    \ell&\mapsto&
                  \displaystyle G(z;z_0)(f\star\rho_{M,\ve})(\ell)
                  =G(z;z_0)\sum_{\substack{n+m= \ell}}f(n)\rho_{M,\ve}(m)
  \end{array}
\end{equation}
This definition is taylored for the next two lemmas.
\begin{lem}
  \label{Approximant1}
  We have $f^\ast(\ell)=0$ when $(\ell,M)\neq1$ and
  $\displaystyle 0\le f^\ast(\ell)\le
  1+\frac{8}{\ve^{|\Xi|}z_0}$.
\end{lem}

\begin{proof}
  Let us notice that any decomposition $\ell=n+m$ with
  $\rho_{M,\ve}(m)\neq0$ has $(n,M)=(\ell,M)$, hence the first
  property.

  By using Section~\ref{ES} and since $\rho_{M,\ve}(m)\ge 0$, we may write
  \begin{equation*}
    f^\ast(\ell)
    \le G(z;z_0)\sum_{\substack{n+m =\ell}}\beta_{z_0,z}(n)\rho_{M,\ve}(m).
  \end{equation*}
  We continue with Lemma~\ref{Fourierbetan}:
  \begin{align*}
    \sum_{\substack{n+m =\ell}}\beta_{z_0,z}(n)\rho_{M,\ve}(m)
    &= \sum_{\substack{q\le z^2}} w_q(z;z_0)
    \sum_{a\mode q}\sum_{m}\rho_{M,\ve}(m)e((\ell-m)a/q)
    \\&= \sum_{\substack{q\le z^2}} \frac{w_q(z;z_0)}{|\Bcal_M(\ve)|^2}
    \sum_{a\mode q}e(\ell a/q)\biggl|\sum_{b\in\Bcal_M(\ve)}e(ba/q)\biggr|^2
    \\&\mkern-50mu= \frac{1}{G(z;z_0)}
    +\sum_{\substack{z_0\le q\le z^2}} \frac{w_q(z;z_0)}{|\Bcal_M(\ve)|^2}
    \sum_{a\mode q}e(\ell a/q)\biggl|\sum_{b\in\Bcal_M(\ve)}e(ba/q)\biggr|^2.
  \end{align*}
  Theorem~\ref{wqLS} applies and gives us the bound
  \begin{equation*}
    f^\ast(\ell)
    \le 1+4\frac{N(Mz_0)^{-1}+z^{2}\log(3z_0)^{-1}}{|\Bcal_M(\ve)|}
    \le 1+\frac{6}{\ve^{|\Xi|}z_0}
    + \frac{6M z^{2}}{N\ve^{|\Xi|}\log(3z_0)}.
  \end{equation*}
  The lemma is proved.
\end{proof}

\begin{lem}
  \label{Approximant2}
  We have
  $\displaystyle
    S(f^\ast,\alpha)
    = G(z;z_0)T^*(\alpha)\bigl|S_M(\alpha,\Xi,\ve)/|\Bcal_M(\ve)|\bigr|^2
  $. 

If, given $\alpha\in\RsurZ$, there exists $\myxi\in\Xi$
  such that $|\alpha-\myxi|\le \ve /N$, 
  then $S(f^\ast,\alpha) =G(z;z_0)T^*(\alpha)(1+\Ocal^*(120\ve))$.

  For any integer $a$, we have $S(f^\ast,a/M) = G(z;z_0)T^*(a/M)$.
\end{lem}

\begin{proof}
  By definition, we have
  \begin{align*}
    S(f^\ast,\alpha)
    &=
    G(z;z_0)\sum_{n,m}f(n)\rho_{M,\ve}(m)e((n+m)\alpha)
   \\&= G(z;z_0)T^*(\alpha)\biggl|\frac{S_M(\alpha,\Xi,\ve)}{|\Bcal_M(\ve)|}\biggr|^2
 \end{align*}
 as claimed. On using Lemma~\ref{GoodPol}, the reader will conclude
 easily, recalling that $\ve\le 1/2$.
\end{proof}
\section{Decomposition}
We have reached the main technical point.
\begin{thm}
  \label{MainDec}
  We have $f=f^\ast G(z;z_0)^{-1}+f^\sharp$ where
  \begin{itemize}
  \item For every $\alpha\in\RsurZ$, $|S(f^\sharp,\alpha)|\le |T^*(\alpha)|$
  and $|S(f^\ast,\alpha)|\le |T^*(\alpha)|G(z;z_0)$.
  \item For every $\alpha\in\RsurZ$, we have $|S(f^\sharp,\alpha)|<
    T^*(0)/A$.
  \item $0\le f^\ast(\ell)\le 1+\epsilon$
    and $f^\ast(\ell)=0$ when $(\ell,M)\neq1$.
  \item For any integer $a$, we have $S(f^\ast,a/M) = G(z;z_0)T^*(a/M)$.
  \end{itemize}
\end{thm}

\begin{proof}
  We have $\ve^{-|\Xi|}\le \epsilon z_0/8$. By construction of $\Xi$,
  every point of $\Cscr(\mathcal{P}^*,A)$ is at distance $\le \ve/N$
  of a point of~$\Xi$. By Lemma~\ref{Approximant2}, this implies that
  $S(f^\sharp,\alpha)=T^*(\alpha)-S(f^\ast,\alpha)G(z;z_0)^{-1}$ satisfies
  \begin{equation*}
    \forall \alpha\in\Cscr(\mathcal{P}^*,A),\quad
    |S(f^\sharp,\alpha)|\le T^*(0)/(2A).
  \end{equation*}
  When $\alpha\notin \Cscr(\mathcal{P}^*,A)$, the inequality
  $|S(f^\sharp,\alpha)|\le |T^*(\alpha)|$ implies that
  $|S(f^\sharp,\alpha)|<T^*(0)/A$ as required.
\end{proof}
\section{Beautification and proof of Theorem~\ref{MainDecB}}
The factor $G(z;z_0)$ in Theorem~\ref{MainDec} may look
awkward to the many. Furthermore, the pointwise bound $f^*\le
1+\epsilon$ does not reflect the fact that the sieve has lost a factor
$(\log N)/G(z;z_0)$. We thus propose to replace $f^\ast$ by
\begin{equation}
  \label{deffflat}
  f^\flat = f^\ast G(z;z_0)^{-1}V(z_0)\log N.
\end{equation}
Only the local upper bound needs to be refreshed.
By Lemma~\ref{Gzz0}, we find that
\begin{equation*}
  f^\flat\le (1+\epsilon)\frac{\log N}{\log (zz_0)}
  \le 2 (1+\epsilon)^2.
\end{equation*}

\begin{proof}[Proof of Theorem~\ref{MainDecB}]
  We have already proved almost everything except that the condition
  $1.02 z_0<\log z$, where $z$ is defined above, holds true. This is equivalent to
  $\exp(2.04 z_0)Mz_0< N$. We already know that $N\ge
  (Mz_0)^{1+\epsilon^{-1}}\ge (Mz_0)^4$ which is why we have imposed
  $\epsilon\le 1/3$. It is thus enough to prove that $M^3\ge \exp(2.04
  z_0)$, which is implied by~$P(z_0)^3\ge \exp(2.04
  z_0)$. The proof of Lemma~\ref{z0vsz} gives us that it is enough to
  show that $\exp(3\frac45 z_0)\ge \exp(2.04 z_0)$ which is obvious.
\end{proof}

\part{Auxiliaries for explicit computations}

\subsubsection*{On  primes}
Let us start with some estimates due to J.B.~Rosser \&
L.~Schoenfeld in~\cite[Theorem 1, Corollary~2,
Theorem 6-8, Theorem~23]{Rosser-Schoenfeld*62}.
\begin{lem}
  \label{RS}
  We have 
  \begin{align*}
    &\prod_{p\le x}\frac{p}{p-1}\le e^\gamma(\log x)
    \biggl(1+\frac{1}{2\log^2x}\biggr)\quad\text{when $x\ge286$},\\
    &e^\gamma(\log x)
      <\prod_{p\le x}\frac{p}{p-1}\le
      e^\gamma(\log x)+\frac{2e^\gamma}{\sqrt{x}}\quad\text{when $x\le 10^8$}.
  \end{align*}
  Futhermore $\pi(x)=\sum_{p\le x}1\le \frac{x}{\log
    x}(1+\frac{3}{2\log x})$ and $\pi(x)\le \frac{5x}{4\log
    x}$, both valid when $x\ge 114$. Finally, $\pi(x)\ge
  x/(\log x)$ when $x\ge17$.
\end{lem}

\begin{lem}
  \label{getVz0}
  When $z_0\ge2$, we have
  $\displaystyle
    \prod_{p< z_0}\frac{p-1}{p}\ge
    \frac{e^{-\gamma}}{\log(9 z_0/5)}$.
    When $z_0>31$, we have
  $\displaystyle
    \prod_{p< z_0}\frac{p-1}{p}\ge
    \frac{e^{-\gamma}}{\log(1.23\, z_0)}$.
\end{lem}
\noindent
The constant $9/5$ is somewhat forced on us by $z_0=3$.
\begin{proof}
  This follows from direct inspection for $z_0\le 100\,000$ and for
  $z_0\le 10^8$ by Lemma~\ref{RS}.
  Again on using this lemma, we find that
  \begin{equation*}
    e^\gamma\log (1.23z_0)\prod_{p< z_0}\frac{p-1}{p}\ge
    \biggl(1+\frac{\log (1.23)}{\log z_0}\biggr)
    \biggl(1+\frac{1}{2\log^2z_0}\biggr)^{-1}.
  \end{equation*}
  The right-hand side is readily seen to be $>1$ when $y=1/\log
  z_0\le0.05$. The lemma follows readily.
\end{proof}

\begin{lem}
  \label{logz0}
  When $z_0\ge 3$, we have
  $\displaystyle\sum_{p<z_0}\frac{\log p}{p-1}\ge \log z_0-0.6$.
\end{lem}

\begin{proof}
  
\end{proof}
\begin{lem}
  \label{z0vsz}
  When $z_0\ge 35$ and $P(z_0)\le z$, we have $z_0\le\tfrac{5}{4}\log z$.
\end{lem}

\begin{proof}
  We have $\vartheta(x)\ge 4x/5$ when $x\ge 30$ (checked up to $10^5$).
\end{proof}
\subsubsection*{On the first $G$-function}

\begin{lem}
  \label{approxGfunctions}
  When $z\ge 1$, we have
  $
    G(z)=\log z+c_0+\Ocal^*((61/25)/\sqrt{z})
    $ where
    \begin{equation*}
      c_0=\gamma+\sum_{p\ge2}\frac{\log p}{p(p-1)}=1.332\,582\,275\,733\ldots
    \end{equation*}
\end{lem}
  This is  part of \cite[Theorem 3.1]{Ramare*18-9}.

\begin{lem}
  \label{PR}
  When $z\ge2$, we have $G(z^2)\le 2G(z)$.
  Also when $z\ge 10$, we have $1.4709\ge G(z)-\log z\ge 1.2$.
\end{lem}

\begin{proof}
  Indeed, the inequality
  \begin{equation*}
    G(z^2)-2G(z)
    \le
    -c_0+\frac{61}{25}\biggl(\frac{1}{\sqrt{z}}+\frac{1}{z}\biggr)
  \end{equation*}
  proves that $G(z^2)-2G(z)\le 0$ when $z\ge 6.5$. When $z\in[n,n+1)$
  and $n$ is an integer, the inequality $G(z^2)\le 2G(z)$ is
  equivalent to by $G(z^2)\le 2G(n)$ and, for this one to hold
  throughout the interval $[n,n+1)$, we need $G(n^2+2n)\le
  2G(n)$. This is readily checked for $n\ge2$. The next upper bound
  for $G(z)$ is~\cite[Lemma 3.5, (1)]{Ramare*95}.  The lower bound
  follows from Lemma~\ref{approxGfunctions} when $z\ge 700$. We
  complete the proof by a direct inspection. The constant $1.2$ is
  forced by the case $z=29$.
\end{proof}

\subsubsection*{On some arithmetic functions}

\begin{lem}
  \label{SQp}
  When real number $Q\ge1$, we have $\displaystyle\sum_{q\le Q}\mu^2(q)\ge Q/2$.
\end{lem}
\noindent
Be cautious: the lower bound of K. Rogers given in \cite{Rogers*64} is
only valid for \emph{integer} values of~$Q$, though this is not specified.
\begin{proof}
  When $Q\ge 1664$, this is a consequence of \cite[Th\'eor\`eme
  3]{Cohen-Dress*88} which asserts that $\sum_{q\le
    Q}\mu^2(q)=6\pi^{-2}Q+\Ocal^*(0.1333\sqrt{Q})$. A straightforward
  numerical check concludes the proof. 
\end{proof}

\begin{lem}
  \label{bettersum}
  When $z_0\ge 35$ and $P(z_0)\le Q^2$, we have
  \begin{equation*}
    \sum_{\substack{q\le Q\\ (q,P(z_0))=1}}1\le 1.1\frac{Q}{\log(3 z_0)}.
  \end{equation*}
\end{lem}

\begin{proof}
  The arithmetical form of the large sieve in its simplest instance
  gives us
  \begin{equation*}
    \sum_{\substack{q\le Q\\ (q,P(z_0))=1}}1\le \frac{Q+z_0^2}{G(z_0)}.
  \end{equation*}
  We readily check, via Lemma~\ref{z0vsz}, that the condition $P(z_0)\le Q^2$ implies that
  $z_0\le Q/\log Q$. Lemma~\ref{PR} provides the lower bound
  $G(z_0)\ge \log (3z_0)$
\end{proof}

\bibliographystyle{plain}
%

\end{document}